\documentclass{amsart}

\usepackage{amsmath} 
\usepackage{amsthm}
\usepackage{amssymb}
\usepackage{graphicx}

\theoremstyle{plain}
\newtheorem{thm}{Theorem}[section]
\newtheorem{prop}[thm]{Proposition}

\newtheorem{lem}[thm]{Lemma}

\theoremstyle{remark}
\newtheorem{rem}[thm]{Remark}

\theoremstyle{definition}

\numberwithin{equation}{section}
\numberwithin{figure}{section}
\numberwithin{table}{section}

\begin{document}

\renewcommand{\thefootnote}{\fnsymbol{footnote}}
\footnote[0]{2000\textit{ Mathematics Subject Classification}
Primary 91B28, 60G10; Secondary 62P05, 93E20.}


\title[\lambdall]{Optimal long term investment model\\
with memory}
\author[\lambdall]{Akihiko Inoue and Yumiharu Nakano}
\address{Department of Mathematics \\
Faculty of Science \\
Hokkaido University \\
Sapporo 060-0810 \\
Japan}
\email{inoue@math.sci.hokudai.ac.jp}
\address{Division of Mathematical Sciences for Social Systems \\
Graduate School of Engineering Science \\ 
Osaka University \\ 
Toyonaka 560-8531, Japan}
\email{y-nakano@sigmath.es.osaka-u.ac.jp}

\date{May 5, 2006}

\keywords{Optimal investment, long term investment, processes with memory, 
processes with stationary increments, Riccati equations, large deviations}

\begin{abstract}
We consider a financial market model driven by an $\mathbf{R}^n$-valued 
Gaussian process with stationary increments which is different from Brownian 
motion. This driving noise process consists of $n$ independent components, and 
each component has memory described by two parameters. 
For this market model, we explicitly solve optimal investment problems. 
These include (i) Merton's portfolio optimization problem; 
(ii) the maximization of growth rate of expected utility of wealth over 
the infinite horizon; 
(iii) the maximization of the large deviation probability 
that the wealth grows at a higher rate than a given benchmark. 
The estimation of paremeters is also considered.
\end{abstract}

\maketitle

\section{Introduction}\label{sec:1}

In this paper we study optimal investment problems for a financial market 
model with memory. This market model $\mathcal{M}$ consists of $n$ risky 
and one riskless assets. The price of the riskless 
asset is denoted by $S_0(t)$ and that of the $i$th risky asset by $S_i(t)$.
We put $S(t)=(S_1(t),\dots,S_n(t))'$, where $A'$ denotes the transpose of a 
matrix $A$. The dynamics of the $\mathbf{R}^n$-valued process 
$S(t)$ are described by the stochastic differential 
equation
\begin{equation}
\begin{split}
dS_i(t)&=S_i(t)\left[\mu_i(t) dt
+\sum\nolimits_{j=1}^n\sigma_{ij}(t)dY_j(t)\right],
\quad t\ge 0, \\
\quad S_i(0)&=s_i,\qquad\qquad\qquad\qquad\qquad\qquad\qquad i=1,\dots,n,
\end{split}
\label{eq:1.1}
\end{equation}
while those of $S_0(t)$ by the ordinary differential equation
\begin{equation}
dS_0(t)=r(t)S_0(t)dt,\quad t\ge 0,\quad S_0(0)=1,
\label{eq:1.2}
\end{equation}
where the coefficients $r(t)\ge 0$, $\mu_i(t)$, and 
$\sigma_{ij}(t)$ are continuous deterministic functions on 
$[0,\infty)$ and the initial prices $s_i$ are positive constants. 
We assume that the $n\times n$ volatility matrix 
$\sigma(t)=(\sigma_{ij}(t))_{1\le i,j\le n}$ is nonsingular for $t\ge 0$. 

The major feature of the model $\mathcal{M}$ is the $\mathbf{R}^n$-valued 
driving noise process $Y(t)=(Y_1(t),\dots,Y_n(t))'$ which has memory. 
We define the $j$th component $Y_j(t)$ by the autoregressive type equation
\begin{equation}
\frac{dY_j(t)}{dt}=-\int_{-\infty }^{t}p_je^{-q_j(t-s)}\frac{dY_j(s)}{ds}ds
+\frac{dW_j(t)}{dt},\quad t\in\mathbf{R},\quad Y_j(0)=0,
\label{eq:1.3}
\end{equation}
where $W(t)=(W_1(t),\dots,W_n(t))'$, $t\in\mathbf{R}$, is an 
$\mathbf{R}^n$-valued standard Brownian motion defined on a complete 
probability 
space $(\Omega,\mathcal{F},P)$, the derivatives $dY_j(t)/dt$ and $dW_j(t)/dt$ 
are in the random distribution sense, and $p_j$'s and $q_j$'s are 
constants such that
\begin{equation}
0<q_j<\infty,\quad -q_j<p_j<\infty,\quad j=1,\dots,n
\label{eq:1.4}
\end{equation}
(cf.\ Anh and Inoue \cite{AI}). 
Equivalently, we may define $Y_j(t)$ by the moving-average type 
representation
\begin{equation}
Y_j(t)=W_j(t)-\int_{0}^{t}\left[ \int_{-\infty }^{s}p_je^{-(q_j+p_j)(s-u)}
dW_j(u)\right]
ds,\quad t\in\mathbf{R}
\label{eq:1.5}
\end{equation}
(see \cite[Examples 2.12 and 2.14]{AI}). 
The components $Y_j(t)$, $j=1,\dots,n$, are Gaussian processes with 
stationary increments 
that are independent of each other. Each $Y_j(t)$ has short 
memory that is described by the two parameters $p_j$ and $q_j$. 
In the special case $p_j=0$, 
$Y_j(t)$ reduces to the Brownian motion $W_j(t)$. 
Driving noise processes with short or long memory of this kind 
are considered in \cite{AI}, Anh et al.\ \cite{AIK} and 
Inoue et al.\ \cite{INA}, for the case $n=1$. 

We define
\[
\mathcal{F}_t:=\sigma\left(\sigma(Y(s): 0\le s\le t)\cup\mathcal{N}\right),
\quad t\ge 0,
\]
where $\mathcal{N}$ is the $P$-null subsets of $\mathcal{F}$. 
This filtration $(\mathcal{F}_t)_{t\ge 0}$ is the underlying 
information structure of the market model $\mathcal{M}$. 
From (\ref{eq:1.5}), we can easily show that 
$(Y(t))_{t\ge 0}$ is a semimartingale with respect to 
$(\mathcal{F}_t)$ (cf.\ \cite[Section 3]{AI}). 
In particular, we can interpret the stochastic differential equation 
(\ref{eq:1.1}) in the usual sense. In actual calculations, however, we need 
explicit semimartingale representations of $Y(t)$. 
It should be noticed that (\ref{eq:1.5}) is not a semimartingale 
representation of $Y(t)$ (except in the special case $p_j=0$). 
For, $W_j(t)$ involves the information of $Y_j(s)$ with $s<0$ 
and vice versa. 
The following two kinds of semimartingale representations of 
$Y(t)$ are obtained in \cite[Example 5.3]{AIK} and \cite[Theorem 2.1]{INA}, 
respectively:
\begin{align}
 Y_j(t)&=B_j(t)-\int_0^t\left[\int_0^s k_j(s,u)dY_j(u)\right]ds,
 \quad t\ge 0,\quad j=1,\dots,n,
\label{eq:1.6} \\
 Y_j(t)&=B_j(t)-\int_0^t\left[\int_0^s l_j(s,u)dB_j(u)\right]ds,
 \quad t\ge 0,\quad j=1,\dots,n,
\label{eq:1.7}
\end{align}
where, for $j=1,\dots,n$, $(B_j(t))_{t\ge 0}$ is 
the so-called \textit{innovation process\/}, i.e., an 
$\mathbf{R}$-valued standard Brownian motion such that
$$
\sigma(Y_j(s): 0\le s\le t)=\sigma(B_j(s): 0\le s\le t),
\quad t\ge 0.
$$
Notice that $B_j$'s are independent of each other. 
The point of \eqref{eq:1.6} and \eqref{eq:1.7} is that 
the deterministic kernels $k_j(t,s)$ and $l_j(t,s)$ are 
given explicitly by
\begin{align}
&k_j(t,s)=p_j(2q_j+p_j)\frac{(2q_j+p_j)e^{q_js}
-p_je^{-q_js}}{(2q_j+p_j)^2e^{q_jt}-p_j^2e^{-q_jt}},
\quad 0\le s\le t,
\label{eq:1.8}\\
&l_j(t,s)=e^{-(p_j+q_j)(t-s)}l_j(s), \quad 0\le s\le t,
\end{align}
with
\begin{equation}
 l_j(s):=p_j\left[1-\frac{2p_jq_j}{(2q_j+p_j)^2e^{2q_js}-p_j^2}\right],
\quad s\ge 0.
\label{eq:1.10}
\end{equation}
We have the equalities
\begin{equation}
\label{eq:1.11}
\int_0^tk_j(t,s)dY_j(s)=\int_0^tl_j(t,s)dB_j(s), \quad t\ge 0, 
\quad j=1,\dots,n.  
\end{equation}

Many authors consider financial market models in which the standard driving 
noise, that is, Brownian motion, is replaced by a different one, such as 
fractional Brownian motion, so that the model can capture 
\textit{memory effect}. 
To name some related contributions, let us mention here 
Comte and Renault \cite{CR1,CR2}, Rogers \cite{R}, Heyde \cite{He}, 
Willinger et al.\ \cite{WTT}, Barndorff-Nielsen and Shephard \cite{BaS}, 
Barndorff-Nielsen et al.\ \cite{BaNS}, Hu and {\O}ksendal \cite{HO}, 
Hu et al.\ \cite{HOS}, Elliott and van der Hoek \cite{EV}, and 
Heyde and Leonenko \cite{HeL}. 
In most of these references, driving noise processes are assumed to have 
stationary increments since this is a natural requirement of 
simplicity. 
Among such models, the above model $\mathcal{M}$ driven by $Y(t)$ which is 
a Gaussian process with {\it stationary increments\/} is possibly 
the simplest one. One advantage of $\mathcal{M}$ is that, by the 
semimartingale representations (\ref{eq:1.6}) and (\ref{eq:1.7}) of $Y(t)$, 
it admits {\it explicit calculations\/} in problems such as those 
considered in 
this paper. Another advantageous feature of the model $\mathcal{M}$ is 
that, assuming $\sigma_{ij}(t)=\sigma_{ij}$, 
real constants, we can easily estimate the characteristic parameters 
$p_j$, $q_j$ and $\sigma_{ij}$ from stock price data. 
We consider this parameter estimation in Appendix C.

For the market model $\mathcal{M}$, we consider an agent who 
has initial endowment $x\in (0,\infty)$ and invests 
$\pi_i(t)X^{x,\pi}(t)$ dollars in the $i$th risky asset for 
$i=1,\dots,n$ and $[1-\sum_{i=1}^{n}\pi_i(t)]X^{x,\pi}(t)$ dollars in 
the riskless asset at each time $t$, where $X^{x,\pi}(t)$ denotes the agent's 
wealth at time $t$. The wealth process $X^{x,\pi}(t)$ is governed by the 
stochastic differential equation
\begin{equation}
\frac{dX^{x,\pi}(t)}{X^{x,\pi}(t)}
=\left[1-\sum\nolimits_{i=1}^n\pi_i(t)\right]\frac{dS_0(t)}{S_0(t)}
 + \sum\nolimits_{i=1}^{n}\pi_i(t)\frac{dS_i(t)}{S_i(t)}, \quad 
X^{x,\pi}(0)=x.
\label{eq:1.12}
\end{equation}
Here, we choose the self-financing strategy 
$\pi(t)=(\pi_1(t),\dots,\pi_n(t))'$ from the admissible class
\[
 \mathcal{A}_T :=\left\{\pi=(\pi(t))_{0\le t\le T}: 
\begin{split}
&\mbox{$\pi$ is an $\mathbf{R}^n$-valued, progressively measurable} \\
&\mbox{process satisfying $\int_{0}^{T}\Vert\pi(t)\Vert^2dt<\infty$ a.s.}
\end{split}
\right\}
\]
for the finite time horizon of length $T\in (0,\infty)$, 
where $\Vert\cdot\Vert$ denotes the Euclidean norm of $\mathbf{R}^n$. 
If the time horizon is infinite, we choose $\pi(t)$ from the class
\[
\mathcal{A}:=\left\{(\pi(t))_{t\ge 0}: (\pi(t))_{0\le t\le T}
 \in\mathcal{A}_T\mbox{ for every }T\in (0,\infty)\right\}. 
\]
Let $\alpha\in (-\infty,1)\setminus\{0\}$ and $c\in\mathbf{R}$. 
In this paper, we consider the following three optimal investment 
problems for the model $\mathcal{M}$:
\begin{align}
&V(T,\alpha):=\sup_{\pi\in\mathcal{A}_T}
  \frac{1}{\alpha}E\left[(X^{x,\pi}(T))^{\alpha}\right],
\tag{\textbf{P1}}\\
&J(\alpha):=\sup_{\pi\in\mathcal{A}} \limsup_{T\to\infty}
  \frac{1}{\alpha T}\log E\left[(X^{x,\pi}(T))^{\alpha}\right],
\tag{\textbf{P2}}\\
&I(c):=\sup_{\pi\in\mathcal{A}} \limsup_{T\to\infty}\frac{1}{T}
\log P\left[X^{x,\pi}(T)\ge e^{cT}\right].
\tag{\textbf{P3}}
\end{align}

The goal of Problem P1 is to maximize the expected utility of 
wealth at the end of the finite horizon. 
This classical optimal investment problem dates back to Merton \cite{Me}. 
We refer to Karatzas and 
Shreve \cite{KS} and references therein for work on this and related problems. 
In Hu et al.\ \cite{HOS}, this problem is solved 
for a Black--Scholes type model driven by fractional Brownian motion. 
In Section \ref{sec:2}, assuming $p_j\ge 0$ for $j=1,\dots,n$, 
we explicitly solve this problem for the model $\mathcal{M}$. 
Our approach is based on a Cameron--Martin 
type formula which we prove in Appendix A. This formula holds under the 
assumption that a relevant Riccati type equation has a solution, and the key 
step of our arguments is to show the existence of such a solution 
(Lemma \ref{lem:2.1}).

The aim of Problem P2 is to maximize the growth rate of 
expected utility of wealth over the infinite horizon. 
This problem is studied by Bielecki and Pliska 
\cite{BP}, and subsequently by other authors under various settings, 
including Fleming and Sheu \cite{FS1,FS2}, Kuroda and Nagai \cite{KN}, 
Pham \cite{P1,P2}, Nagai and Peng \cite{NP}, Hata and Iida \cite{HI}, 
and Hata and Sekine \cite{HS1,HS2}. 
In Section \ref{sec:3}, we solve Problem P2 for the model $\mathcal{M}$ by 
verifying that a 
candidate of optimal strategy suggested by the solution to Problem P1 
is actually optimal. In so doing, 
existence results on solutions to 
Riccati type equations (Lemmas \ref{lem:2.1} and \ref{lem:3.5}) 
play a key role as in Problem P1. 
The result of Nagai and Peng \cite{NP} on the 
asymptotic behavior of solutions to Riccati equations, 
which we review in Appendix B, is also an essential ingredient in 
our arguments.

The purpose of Problem P3 is to maximize the large deviation probability that 
the wealth grows at a higher rate than 
the given benchmark $c$. 
This problem is studied by Pham 
\cite{P1,P2}, in which a significant result, that is, a duality relation 
between Problems P2 and P3, is established. 
Subsequently, this problem is studied by Hata and Iida 
\cite{HI} and Hata and Sekine \cite{HS1,HS2} under different settings. 
In Section \ref{sec:4}, we solve Problem P3 for the market 
model $\mathcal{M}$. In the approach of \cite{P1,P2}, 
one needs an explicit expression of $J(\alpha)$. 
Since our solution to Problem P2 is explicit, we can solve Problem 
P3 for $\mathcal{M}$ using this approach. 
As in \cite{P1,P2}, our solution to Problem 3 is given in the form of a 
sequence of nearly optimal strategies. 
For $c<\bar{c}$ with certain constant $\bar{c}$, 
an optimal strategy, rather than such a nearly optimal sequence, 
is obtained by ergodic arguments.


\section{Optimal investment over the finite horizon}\label{sec:2}

In this section, we consider the finite horizon optimization 
problem P1 for the market model $\mathcal{M}$. 
Throughout this section, we assume 
$\alpha\in (-\infty,1)\setminus\{0\}$ and
\begin{equation}
0<q_j<\infty,\quad 0\le p_j<\infty,\quad j=1,\dots,n.
\label{eq:2.1}
\end{equation}
Thus $p_j\ge 0$ rather than $p_j>-q_j$ (see Remark \ref{rem:2.6} below).

Let $Y(t)=(Y_1(t),\dots,Y_n(t))'$ and $B(t)=(B_1(t),\dots,B_n(t))'$ be 
the driving noise and innovation processes, 
respectively, described in Section \ref{sec:1}. 
We define an $\mathbf{R}^n$-valued deterministic function 
$\lambda(t)=(\lambda_1(t),\dots,\lambda_n(t))'$ by
\begin{equation}
\lambda(t)
:=\sigma^{-1}(t)\left[\mu(t)-r(t)\mathbf{1}\right],
\quad t\ge 0,
\label{eq:2.2}
\end{equation}
where $\mathbf{1}:=(1,\dots,1)'\in\mathbf{R}^n$. 
For the kernels $k_j(t,s)$'s in \eqref{eq:1.8}, we put 
\[
k(t,s):=\mathrm{diag}(k_1(t,s),\dots,k_n(t,s)), \quad 0\le s\le t.
\]
We denote by $\xi(t)=(\xi_1(t),\dots,\xi_n(t))'$ the $\mathbf{R}^{n}$-valued 
process $\int_0^tk(t,s)dY(s)$, 
i.e.,
\begin{equation}
\label{eq:2.3}
\xi_j(t):=\int_0^tk_j(t,s)dY_j(s),  \quad t\ge 0,\quad j=1,\dots,n.  
\end{equation}
By \eqref{eq:1.1}, \eqref{eq:1.2}, \eqref{eq:1.6}, and \eqref{eq:1.12}, 
the wealth process $X^{x,\pi}(t)$ evolves according to 
\[
\frac{dX^{x,\pi}(t)}{X^{x,\pi}(t)}
=r(t)dt+\pi'(t)\sigma(t)\left[\lambda(t)-\xi(t)\right]dt
 +\pi'(t)\sigma(t)dB(t),\quad t\ge 0,
\]
whence, by the It\^{o} formula, we have, for $t\ge 0$,
\begin{equation}
\begin{split}
X^{x,\pi}(t)
&=x\exp\left[\int_0^t \left\{r(s) + \pi'(s)\sigma(s)
\left(\lambda(s)-\xi(s)\right)- 
\frac{1}{2}\Vert\sigma'(s)\pi(s)\Vert^2 \right\}ds\right. \\
&\qquad\qquad\qquad\qquad\qquad\qquad\qquad
+ \left.\int_0^t \pi'(s)\sigma(s)dB(s)\right].
\end{split}
\label{eq:2.4}
\end{equation}

We define an $\mathbf{R}$-valued process $Z(t)$ by
\[
Z(t):=\exp\left[-\int_0^t\left\{\lambda(s)-\xi(s)\right\}'dB(s)
-\frac{1}{2}\int_0^t\left\Vert\lambda(s)
-\xi(s)\right\Vert^2ds\right], \quad t\ge 0.
\]
Since $\lambda(t)-\xi(t)$ is a continuous Gaussian 
process, the process $Z(t)$ is a $P$-martingale 
(see, e.g., Example 3(a) in Liptser and Shiryayev \cite[Section 6.2]{LS}). 
We define the $\mathbf{R}$-valued process $(\Gamma(t))_{0\le t\le T}$ by
\[
\Gamma(t):=E\left[\left.Z^{\beta}(T)\right|\mathcal{F}_t\right],
\quad 0\le t\le T,
\]
where $\beta$ is the 
conjugate exponent of $\alpha$, i.e.,
\[
(1/\alpha)+(1/\beta)=1.
\]
Notice that $0<\beta<1$ (resp.\ $-\infty<\beta<0$) if 
$-\infty<\alpha<0$ (resp.\ $0<\alpha<1$). 
In view of Theorem 7.6 in Karatzas and Shreve \cite[Chapter 3]{KS}, 
to solve Problem P1, we only have to derive a stochastic integral 
representation for $\Gamma(t)$.

We define an $\mathbf{R}$-valued $P$-martingale $K(t)$ by
\[
K(t):=\exp\left[-\beta\int_0^t \{\lambda(s)-\xi(s)\}'dB(s)
-\frac{\beta^2}{2}\int_0^t \Vert\lambda(s)-\xi(s)\Vert^2ds\right],
\quad t\ge 0.
\]
Then, by Bayes' rule, we have
\[
\begin{split}
\Gamma(t)&=E\left[\left. K(T)\exp\left\{-\frac{1}{2}\beta(1-\beta)
 \int_0^T \Vert\lambda(s)-\xi(s)\Vert^2ds\right\}\right|
 \mathcal{F}_t\right] \\
&= K(t)\bar{E}\left[\left.\exp\left\{-\frac{1}{2}\beta(1-\beta)
 \int_0^T \Vert\lambda(s)-\xi(s)\Vert^2ds\right\}\right|
 \mathcal{F}_t\right]
\end{split}
\]
for $t\in [0,T]$, where $\bar{E}$ stands for the expectation with respect to 
the probability measure $\bar{P}$ on $(\Omega,\mathcal{F}_T)$ such that 
$d\bar{P}/dP=K(T)$. 
Thus
\begin{equation}
\begin{split}
\Gamma(t)&= Z^{\beta}(t)\exp\left\{-\frac{1}{2}\beta
 (1-\beta)\int_t^T \Vert\lambda(s)\Vert^2ds\right\}\\
&\quad\times \bar{E}\left[\left.\exp\left\{-\frac{1}{2}\beta
 (1-\beta)\int_t^T \left(\Vert\xi(s)\Vert^2-2\lambda'(s)\xi(s)\right)ds
\right\}\right|\mathcal{F}_t\right].
\end{split}
\label{eq:2.5}
\end{equation}

We are to apply Theorem \ref{thm:A1} in Appendix A to \eqref{eq:2.5}. 
By \eqref{eq:1.11}, the dynamics of $\xi(t)$ are described by 
the $n$-dimensional stochastic differential equation
\begin{equation}
d\xi(t)=-(p+q)\xi(t)dt+l(t)dB(t),\quad t\ge 0,
\label{eq:2.6}
\end{equation}
where $p:=\mathrm{diag}(p_1,\dots,p_n)$, 
$q:=\mathrm{diag}(q_1,\dots,q_n)$, and 
$l(t):=\mathrm{diag}(l_1(t),\dots,l_n(t))$ 
with $l_j(t)$'s as in (\ref{eq:1.10}). 
Write $\bar{B}(t):=B(t)+\beta\int_0^t[\lambda(s)-\xi(s)]ds$ for 
$t\in [0, T]$. Then $\bar{B}(t)$ is an $\mathbf{R}^n$-valued standard Brownian 
motion under $\bar{P}$. By \eqref{eq:2.6}, the process $\xi(t)$ evolves 
according to
\begin{equation}
d\xi(t)=\left[\rho(t)+b(t)\xi(t)\right]dt+l(t)d\bar{B}(t),
\quad t\ge 0,
\label{eq:2.7}
\end{equation}
where $\rho(t)=(\rho_1(t),\dots,\rho_n(t))'$, 
$b(t)=\mathrm{diag}(b_1(t),\dots,b_n(t))$ with
\begin{align}
&\rho_j(t):=-\beta l_j(t)\lambda_j(t), \quad\quad t\ge 0,\quad j=1,\dots,n,
\label{eq:2.8}\\
&b_j(t):=-(p_j+q_j)+\beta l_j(t), \quad t\ge 0,\quad j=1,\dots,n.
\label{eq:2.9}
\end{align}
By Theorem \ref{thm:A1} in Appendix A, we are led to consider 
the following one-dimensional backward Riccati equations: 
for $j=1,\dots,n$
\begin{equation}
\dot{R}_j(t)-l_j^2(t)R_j^2(t)+2b_j(t)R_j(t)+\beta(1-\beta)=0,
\quad 0\le t\le T,\quad R_j(T)=0.
\label{eq:2.10}
\end{equation}

The following lemma, especially (iii), is crucial in our arguments.

\begin{lem}\label{lem:2.1}Let $j\in\{1,\dots,n\}$.
\begin{enumerate}
\item If $p_j=0$, then {\rm \eqref{eq:2.10}} has a unique 
solution $R_j(t)\equiv R_j(t;T)$.
\item If $-\infty<\alpha<0$, then 
{\rm \eqref{eq:2.10}} has a unique nonnegative solution $R_j(t)\equiv 
R_j(t;T)$.  
\item If $p_j>0$ and $0<\alpha<1$, then 
{\rm (\ref{eq:2.10})} has a unique solution $R_j(t)\equiv R_j(t;T)$ such that 
$R_j(t)\ge b_{j}(t)/l_{j}^2(t)$ for $t\in [0,T]$.
\end{enumerate}
\end{lem}

\begin{proof}(i)\ If $p_j=0$, then (\ref{eq:2.10}) is linear, 
whence it has a unique solution.

(ii)\ If $-\infty<\alpha<0$, then $\beta(1-\beta)>0$, so that, 
by the well-known result on Riccati 
equations (see, e.g., Fleming and Rishel \cite[Theorem 5.2]{FR} and 
Liptser and Shiryayev \cite[Theorem 10.2]{LS}), 
(\ref{eq:2.10}) has a unique nonnegative solution.

(iii)\ When $p_j>0$ and $0<\alpha<1$, write
\begin{equation}
a_{1}(t):=l_j^2(t),\quad 
a_{2}(t):=b_j(t),\quad 
a_3:=\beta(1-\beta),\quad t\ge 0.
\label{eq:2.11}
\end{equation}
Then the equation for $P(t):=R_j(t)-[a_{2}(t)/a_{1}(t)]$ becomes
\begin{equation}
\dot{P}(t)-a_{1}(t)P^2(t)+a_{4}(t)=0,\quad 0\le t\le T,
\label{eq:2.12}
\end{equation}
where
\[
a_{4}(t):=\frac{a_{2}^2(t)+a_{1}(t)a_3}{a_{1}(t)}
+\frac{d}{dt}\left[\frac{a_{2}(t)}{a_{1}(t)}\right].
\]
Since $dl_j(t)/dt>0$ and $\beta<0$, we see that
\[
\frac{d}{dt}\left[\frac{a_{2}(t)}{a_{1}(t)}\right]
=\frac{2(p_j+q_j)-\beta l_j(t)}{l_j(t)^3}\cdot\frac{dl_j}{dt}(t)>0.
\]
We write $a_{2}^2(t)+a_{1}(t)a_3$ as
\[
(1-\beta)\left[(p_j+q_j)^2-\{(p_j+q_j)-l_j(t)\}^2\right]
+[(p_j+q_j)-l_j(t)]^2,
\]
which is positive since $0\le l_j(t)\le p_j$. 
Thus $a_{4}(t)>0$, so that 
(\ref{eq:2.12}) has a unique nonnegative solution $P(t)\equiv P(t;T)$. 
The desired solution to (\ref{eq:2.10}) 
is given by $R_j(t)=P(t)+[a_{2}(t)/a_{1}(t)]$.
\end{proof}

In what follows, we write $R_j(t)\equiv R_j(t;T)$ for the unique solution to 
(\ref{eq:2.10}) in the sense of Lemma \ref{lem:2.1}. 
Then $R(t):=\mathrm{diag}(R_1(t),\dots,R_n(t))$ satisfies 
the backward matrix Riccati equation
\begin{equation}
\begin{split}
&\dot{R}(t)-R(t)l^2(t)R(t)+b(t)R(t)+R(t)b(t)+ \beta(1-\beta)I_n=0,\quad 
0\le t\le T,\\
&R(T)=0,
\end{split}
\label{eq:2.13}
\end{equation}
where $I_n$ denotes the $n\times n$ unit matrix. 
For $j=1,\dots,n$, let $v_j(t)\equiv v_j(t;T)$ be the solution to the 
following one-dimensional linear equation:
\begin{equation}
\begin{split}
&\dot{v}_j(t)+[b_j(t)-l_j^2(t)R_j(t;T)]v_j(t)+
 \beta(1-\beta)\lambda_j(t)-R_j(t;T)\rho_j(t)=0,\\
&\quad0\le t\le T,\quad v_j(T)=0.
\end{split}
\label{eq:2.14}
\end{equation}
Then $v(t)\equiv v(t;T):=(v_1(t;T),\dots,v_n(t;T))'$ satisfies 
the matrix equation
\begin{equation}
\begin{split}
&\dot{v}(t)+[b(t)-l^2(t)R(t;T)]v(t)+
 \beta(1-\beta)\lambda(t)-R(t;T)\rho(t)=0,\\
&\quad0\le t\le T,\quad v(T)=0.
\end{split}
\label{eq:2.15}
\end{equation}
We put, for $j=1,\dots,n$ and $(t,T)\in \Delta$,
\begin{equation}
g_j(t;T):= v_j^2(t;T)l_j^2(t)+2\rho_j(t)v_j(t;T)-l_j^2(t)R_j(t;T)
 -\beta(1-\beta)\lambda_j^2(t),
\label{eq:2.16}
\end{equation}
where
\begin{equation}
\Delta:=\{(t,T): 0<T<\infty,\ 0\le t\le T\}.
\label{eq:2.17}
\end{equation}

We are now ready to give the desired representation for $\Gamma(t)$.

\begin{prop}\label{prop:2.2}
Write
\begin{equation}
\psi(t):=\Gamma(t)\left[-\beta\lambda(t)+
 \{\beta-l(t)R(t;T)\}\xi(t)+l(t)v(t;T)\right],
\quad 0\le t\le T.
\label{eq:2.18}
\end{equation}
Then, for $t\in [0,T]$, we have 
$\Gamma(t)=\Gamma(0)+\int_0^t\psi'(s)dB(s)$ with
\begin{equation}
\Gamma(0)
=\exp\left[\frac{1}{2}\int_0^T \sum\nolimits_{j=1}^{n} g_j(s;T)ds\right].
\label{eq:2.19}
\end{equation} 
\end{prop}

\begin{proof}
It follows from \eqref{eq:2.5}, \eqref{eq:2.7}, \eqref{eq:2.13}, 
\eqref{eq:2.15} and Theorem \ref{thm:A1} that
\begin{equation}
\Gamma(t)=Z^{\beta}(t)\exp\left[\sum_{j=1}^n
\left\{v_j(t)\xi_j(t)-\frac{1}{2}\xi_j^2(t)R_j(t)
  +\frac{1}{2}\int_t^T g_j(s;T)ds\right\}\right]. 
\label{eq:2.20}
\end{equation}
The equality (\ref{eq:2.19}) follows from this. 
A straightforward calculation based on \eqref{eq:2.20}, \eqref{eq:2.6} and 
the It\^{o} formula gives $d\Gamma(t)=\psi'(t)dB(t)$, 
where $\psi(t)$ is as in \eqref{eq:2.18}. 
Thus the proposition follows. 
\end{proof}

Recall that we have assumed $\alpha\in (-\infty,1)\setminus\{0\}$ and 
\eqref{eq:2.14}. Here is the solution to Problem P1.

\begin{thm}
\label{thm:2.3}
For $T\in (0,\infty)$, the strategy 
$(\hat{\pi}_T(t))_{0\le t\le T}\in\mathcal{A}_T$ defined by
\begin{equation}
\hat{\pi}_T(t):=(\sigma')^{-1}(t)
 \left[(1-\beta)\{\lambda(t)-\xi(t)\}-l(t)R(t;T)\xi(t)
 +l(t)v(t;T)\right]
\label{eq:2.21}
\end{equation}
is the unique optimal strategy for Problem P1. 
The value function $V(T)\equiv V(T,\alpha)$ in {\rm (P1)} is given by
\begin{equation}
V(T)=\frac{1}{\alpha}[xS_0(T)]^{\alpha}
\exp\left[\frac{(1-\alpha)}{2}
\sum\nolimits_{j=1}^{n}\int_0^T g_j(t;T)dt\right].
\label{eq:2.22}
\end{equation}
\end{thm}

\begin{proof}
By Theorem 7.6 in Karatzas and Shreve \cite[Chapter 3]{KS}, 
the unique optimal strategy $\pi_T(t)$ for Problem P1 is given by
\[
\pi_T(t):=(\sigma')^{-1}(t)
\left[\Gamma^{-1}(t)\psi(t)+\lambda(t)-\xi(t)\right],
\quad 0\le t\le T, 
\]
which, by (\ref{eq:2.18}), is equal to $\hat{\pi}_T(t)$. Thus 
the first assertion follows. By the same theorem in 
\cite{KS}, 
$V(T)=\alpha^{-1}[xS_0(T)]^{\alpha}\Gamma^{1-\alpha}(0)$. 
This and (\ref{eq:2.19}) give (\ref{eq:2.22}).
\end{proof}

\begin{rem}
\label{rem:2.4}
We can regard $\xi(t)=\int_0^tk(t,s)dY(s)$, which is the only 
random term on the right-hand side of \eqref{eq:2.21}, 
as representing the memory effect. 
To illustrate this point, suppose that $(\sigma_{ij}(t))$ is 
a constant matrix. 
Then, by \eqref{eq:C2} in Appendix C, we can express $Y(t)$, whence $\xi(t)$, 
in terms of the past prices $S(u)$, $u\in [0,t]$, of the risky assets.
\end{rem}

\begin{rem}
\label{rem:2.5}
From \cite[Theorem 7.6]{KS}, we also find that
\[
X^{x,\hat{\pi}_T}(t)=x\frac{S_0(t)\Gamma(t)}{Z(t)\Gamma(0)},
\quad 0\le t\le T.
\]   
\end{rem}

\begin{rem}
\label{rem:2.6}
Regarding \eqref{eq:2.1}, we assume this to ensure the existence of 
solution to \eqref{eq:2.10} for $j=1,\dots,n$. 
Under the weaker assumption \eqref{eq:1.4}, 
we could show by a different argument that, for $j=1,\dots,n$, 
\eqref{eq:2.10} has a solution 
if $\alpha\in(-\infty, \bar{\alpha}_j)\setminus\{0\}$, where 
$\bar{\alpha}_j\in (0,1]$ is defined by
\[
\bar{\alpha}_j:=1\quad \mbox{if} \ 0\le p_j<\infty,\quad 
:=\frac{(p_j+q_j)^2}{l_j^2(0)+q_j^2}\quad \mbox{if} \ -q_j<p_j<0.
\]
From this, we see that 
the same result as Theorem \ref{thm:2.3} holds under \eqref{eq:1.4} if 
$-\infty<\alpha<\bar{\alpha}$, $\alpha\ne 0$, where 
$\bar{\alpha}:=\min\{\bar{\alpha}_j: j=1,\dots,n\}$. 
However, we did not succeed in extending the result to the 
most general case $-\infty<\alpha<1$, $\alpha\ne 0$. 
Such an extension, if possible, would lead us to the solution of 
Problem P3 under \eqref{eq:1.4} (see Remark \ref{rem:3.8}).
\end{rem}


\section{Optimal investment over the infinite horizon}\label{sec:3}

In this section, we consider the infinite horizon optimization problem P2 
for the financial market model $\mathcal{M}$. 
Throughout this section, we assume \eqref{eq:2.1} and the following two 
conditions:
\begin{align}
&\lim_{T\to\infty}\frac{1}{T}\int_0^T r(t)dt=\bar{r}\quad 
\mbox{with}\;\; \bar{r}\in [0,\infty),
\label{eq:3.1}\\
&\lim_{t\to\infty}\lambda(t)=\bar{\lambda}\quad 
\mbox{with}\;\; \bar{\lambda}=(\bar{\lambda}_1,\dots,\bar{\lambda}_n)'\in 
\mathbf{R}^n.
\label{eq:3.2}
\end{align}
Here recall $\lambda(t)=(\lambda_1(t),\dots,\lambda_n(t))'$ from 
(\ref{eq:2.2}). In the main result of this section (Theorem \ref{thm:3.4}), 
we will also assume 
$\alpha^*<\alpha<1$, $\alpha\ne 0$, where
\begin{equation}
\alpha^*:=\max(\alpha_1^*,\dots,\alpha_n^*)
\label{eq:3.3}
\end{equation}
with
\begin{equation}
\alpha^*_j:=\left\{
\begin{split}
&-\infty \quad\mbox{if}\quad 0\le p_j\le 2q_j, \\
&-3-\frac{8q_j}{p_j-2q_j} \quad\mbox{if}\quad 2q_j<p_j<\infty.
\end{split}
\right.
\label{eq:3.4}
\end{equation}
Notice that $\alpha^*\in [-\infty,-3)$.

To give the solution to Problem P2, we take the following steps:
\begin{enumerate}
\item For the value function $V(T)\equiv V(T,\alpha)$ in (P1), we calculate 
the following limit explicitly:
\begin{equation}
\tilde{J}(\alpha):=\lim_{T\to\infty}\frac{1}{\alpha T}\log [\alpha V(T)].
\label{eq:3.5}
\end{equation}
\item For $\hat{\pi}\in\mathcal{A}$ in \eqref{eq:3.14} below, 
we calculate the growth rate
\begin{equation}
J^*(\alpha)
:=\lim_{T\to\infty}
\frac{1}{\alpha T}\log E\left[(X^{x,\hat{\pi}}(T))^{\alpha}\right],
\label{eq:3.6}
\end{equation}
and verify that $J^*(\alpha)=\tilde{J}(\alpha)$.
\item Since the definition of $V(T)$ implies
\begin{equation}
\limsup_{T\to\infty}\frac{1}{\alpha T}\log E[(X^{x,\pi}(T))^{\alpha}]
\le \tilde{J}(\alpha),\quad \forall\pi\in\mathcal{A},
\label{eq:3.7}
\end{equation}
we conclude that $\hat{\pi}$ is an optimal strategy for Problem P2 and that 
the optimal growth rate $J(\alpha)$ in (P2) is given by 
$J(\alpha)=J^*(\alpha)=\tilde{J}(\alpha)$.
\end{enumerate}

Let $\alpha\in (-\infty,1)\setminus\{0\}$ and $\beta$ be its conjugate 
exponent as in Section \ref{sec:2}. 
For $j=1,\dots,n$, recall $b_j(t)$ from (\ref{eq:2.9}). 
We have $\lim_{t\to\infty}b_j(t)=\bar{b}_j$, where
\[
 \bar{b}_j:=-(1-\beta)p_j-q_j.
\]
Notice that $\bar{b}_j<0$. We consider the equation
\begin{equation}
p_j^2x^2-2\bar{b}_jx-\beta(1-\beta)=0.
\label{eq:3.8}
\end{equation}
When $p_j=0$, we write $\bar{R}_j$ 
for the unique solution $\beta(1-\beta)/(2q_j)$ of this linear 
equation. If $p_j>0$, then 
\[
\bar{b}_j^2+\beta(1-\beta)p_j^2=(1-\beta)[(p_j+q_j)^2-q_j^2]+q_j^2\ge q_j^2>0,
\]
so that we may write $\bar{R}_j$ for the larger solution to the quadratic 
equation \eqref{eq:3.8}. Let 
\[
K_j:=\sqrt{\bar{b}_j^2+\beta(1-\beta)p_j^2}.
\]
Then $\bar{b}_j-p_j^2\bar{R}_j=-K_j<0$.

As in Section \ref{sec:2}, 
we write $R_j(t)\equiv R_j(t;T)$ for the unique solution to 
(\ref{eq:2.10}) in the sense of Lemma \ref{lem:2.1}. 
Recall $\Delta$ from (\ref{eq:2.17}). 
The next proposition provides the necessary results on the asymptotic 
behavior of $R_j(t;T)$.

\begin{prop}\label{prop:3.1}
Let $-\infty<\alpha<1$, $\alpha\ne 0$, and 
$j\in\{1,\dots,n\}$. Then
\begin{enumerate}
\item $R_j(t;T)$ is bounded in $\Delta$.
\item $\lim_{T-t\to\infty,\ t\to\infty}R_j(t;T)=\bar{R}_j$.
\item For $\delta, \epsilon\in (0,\infty)$ such that $\delta + \epsilon<1$, 
\[
\lim_{T\to\infty}\sup_{\delta T\le t\le (1-\epsilon)T}
 |R_j(t;T)-\bar{R}_j|=0.
\]
\end{enumerate}
\end{prop}

\begin{proof}
If $p_j=0$, then $l_j(t)=0$ and $b_j(t)=-q_j<0$ for $t\ge 0$, 
so that the assertions follow from Theorem \ref{thm:B3} in Appendix B. 

We assume $p_j>0$. 
Since
\[
 \vert l_j(t)-p_j\vert\le \frac{p_j^2}{2(p_j+q_j)}e^{-2q_jt}, \quad t\ge 0,
\]
the function $l_j(t)$ converges to $p_j$ exponentially fast as $t\to\infty$. 
Hence the coefficients of the equation (\ref{eq:2.10}) converge to 
their counterparts in (\ref{eq:3.8}) exponentially fast, too. 
If $-\infty<\alpha<0$, then the desired 
assertions follow from Theorem \ref{thm:B1} in Appendix B 
(due to Nagai and Peng \cite{NP}). 
Suppose $0<\alpha<1$. 
Let $a_1(t)$, $a_2(t)$ and $a_3$ be as in (\ref{eq:2.11}). 
Since $R_j(t;T)\ge b_{j}(t)/l_{j}(t)^2$ and $b_{j}(t)/l_{j}(t)^2$ is bounded 
from below in $\Delta$, so is 
$R_j(t;T)$. To show that 
$R_j(t;T)$ is bounded from above in $\Delta$, 
we consider the solution $M_j(t)\equiv M_j(t;T)$ to the linear equation
\[
 \dot{M}_j(t)+2[a_2(t)-\bar{R}_ja_1(t)]M_j(t)+a_3+a_1(t)\bar{R}^2_j=0,
\quad 0\le t\le T,\quad M_j(T)=0.
\]
Since $M_j(T)-R_j(T)=0$ and
\[
 [\dot{M}_j(t)-\dot{R}_j(t)]+2[a_2(t)-\bar{R}_ja_1(t)][M_j(t)-R_j(t)]
 =-a_1(t)\left[R_j(t)-\bar R_j\right]^2\le 0,
\]
we have $R_j(t;T)\le M_j(t;T)$ in $\Delta$. However, 
$a_2(t)-\bar{R}_ja_1(t)\to 
\bar{b}_j-\bar{R}_j\bar{p}_j^2<0$ as $t\to\infty$, 
so that $M_j(t;T)$ is bounded from above in $\Delta$, whence 
so is $R_j(t;T)$. The desired 
assertions now follow from Theorem \ref{thm:B2} in Appendix B.
\end{proof}

Let $j\in\{1,\dots,n\}$. For $\rho_j(t)$ in (\ref{eq:2.8}), 
we have $\lim_{t\to\infty}\rho_j(t)=\bar{\rho}_j$, where
\[
\bar{\rho}_j:=-\beta p_j\bar{\lambda}_j.
\]
Let $v_j(t)\equiv v_j(t;T)$ be the solution to \eqref{eq:2.14} as in 
Section \ref{sec:2}. Define $\bar{v}_j$ by
\begin{equation}
\left(\bar{b}_j-p_j^2\bar{R}_j\right)\bar{v}_j+\beta(1-\beta)\bar{\lambda}_j
-\bar{R}_j\bar{\rho}_j=0.
\label{eq:3.9}
\end{equation}

\begin{prop}\label{prop:3.2}
Let $-\infty<\alpha<1$, $\alpha\ne 0$, and 
$j\in\{1,\dots,n\}$. Then
\begin{enumerate}
\item $v_j(t;T)$ is bounded in $\Delta$.
\item $\lim_{T-t\to\infty,\ t\to\infty}v_j(t;T)=\bar{v}_j$.
\item For $\delta, \epsilon\in (0,\infty)$ 
such that $\delta + \epsilon<1$, 
\[
\lim_{T\to\infty}\sup_{\delta T\le t\le (1-\epsilon)T}
\vert v_j(t;T)-\bar{v}_j\vert=0.
\]
\end{enumerate}
\end{prop}

\begin{proof}
The coefficients of (\ref{eq:2.14}) converge to their counterparts 
in (\ref{eq:3.9}). Also,
\[
\lim_{T-t\to\infty, \; t\to\infty}[b_j(t)-l_j^2(t)R_j(t;T)]
 = \bar{b}_j-p_j^2\bar{R}_j=-K_j<0.
\]
Thus the proposition follows from Theorem \ref{thm:B3} in Appendix B.
\end{proof}

For $j=1,\dots,n$ and $-\infty<\alpha<1$, $\alpha\ne 0$, we put
\begin{gather}
F_j(\alpha):=\frac{(p_j+q_j)^2\bar{\lambda}_j^2\alpha}
{\left[(1-\alpha)(p_j+q_j)^2+\alpha p_j(p_j+2q_j)\right]},
\label{eq:3.10}\\
\begin{split}
&G_j(\alpha)\\
&\ \ :=(p_j+q_j)-q_j\alpha-(1-\alpha)^{1/2}
\left[(1-\alpha)(p_j+q_j)^2+\alpha p_j(p_j+2q_j)\right]^{1/2}.
\end{split}
\label{eq:3.11}
\end{gather}
Recall the value function $V(T)\equiv V(T,\alpha)$ from (P1) and 
its representation (\ref{eq:2.22}). 
In the next proposition, we compute $\tilde{J}(\alpha)$ in \eqref{eq:3.5}.

\begin{prop}
\label{prop:3.3}
Let $-\infty<\alpha<1$, $\alpha\ne 0$. Then the limit 
$\tilde{J}(\alpha)$ in \eqref{eq:3.5} exists and is given by 
\begin{equation}
\tilde{J}(\alpha)=\bar{r}+\frac{(1-\alpha)}{2\alpha}
\sum_{j=1}^{n}\bar{g}_j,
\label{eq:3.12}
\end{equation}
where
\[
 \bar{g}_j:=\bar{v}^2_jp_j^2+2\bar{\rho}_j\bar{v}_j-p_j^2\bar{R}_j
 -\beta(1-\beta)\bar{\lambda}_j^2,\quad j=1,\dots,n.
\]
More explicitly,
\begin{equation}
\tilde{J}(\alpha)=\bar{r}+\frac{1}{2\alpha}\sum_{j=1}^nF_j(\alpha)
+\frac{1}{2\alpha}\sum_{j=1}^{n}G_j(\alpha).
\label{eq:3.13}
\end{equation}
\end{prop}

\begin{proof}
Recall $g_j(t;T)$ from (\ref{eq:2.16}). 
By Propositions \ref{prop:3.1} and \ref{prop:3.2}, 
\[
\lim_{T\to\infty}\frac{1}{T}\int_0^T g_j(t;T)dt=\bar{g}_j,
\quad j=1,\dots,n.
\]
From this and (\ref{eq:2.22}),
\[
\begin{split}
 \frac{1}{\alpha T}\log\left[\alpha V(T)\right]
 &=\frac{\log x}{T} + \frac{1}{T}\int_0^T r(t)dt
   +\frac{1-\alpha}{2\alpha}\sum_{j=1}^n\frac{1}{T}\int_0^T g_j(t;T)dt\\
&\to \bar{r}+\frac{1-\alpha}{2\alpha}\sum_{j=1}^{n}\bar{g}_j
 \quad \mbox{as $T\to\infty$,}
\end{split}
\]
which implies \eqref{eq:3.12}.

We have $\bar{v}_j=\beta\bar{\lambda}_j(1-\beta+p_j\bar{R}_j)/K_j$. 
Also,
\[
\begin{split}
\beta p_j^2(p_j\bar{R}_j)^2-2\beta p_j^2\bar{R}_jK_j
&=\beta p_j^2\bar{R}_j(p_j^2\bar{R}_j-2K_j)
=\beta p_j^2\bar{R}_j(\bar{b}_j-K_j)\\
&=\beta(\bar{b}_j^2-K_j^2)=-\beta^2(1-\beta)p_j^2,
\end{split}
\]
and $\beta(1-\beta)=-\alpha/(1-\alpha)^2$. Thus
\[
\begin{split}
\bar{v}^2_jp_j^2&+2\bar{\rho}_j\bar{v}_j
 -\beta(1-\beta)\bar{\lambda}_j^2\\
&=\frac{\beta\bar{\lambda}_j^2}{K_j^2}
\left[\beta p_j^2(1-\beta+p_j\bar{R}_j)^2-2\beta p_j(1-\beta+p_j\bar{R}_j)K_j
-(1-\beta)K_j^2\right]\\
&=\frac{\beta\bar{\lambda}_j^2}{K_j^2}
\left[\{\beta p_j^2(p_j\bar{R}_j)^2-2\beta p_j^2\bar{R}_jK_j\}
+2\beta(1-\beta)p_j(p_j^2\bar{R}_j-K_j)\right.\\
&\left.\qquad\qquad\qquad\qquad\qquad\qquad\qquad\qquad\quad 
+\beta p_j^2(1-\beta)^2-(1-\beta)K_j^2\right]\\
&=\frac{\beta(1-\beta)\bar{\lambda}_j^2}{K_j^2}
\left[-\beta^2p_j^2+2\beta p_j\bar{b}_j+\beta(1-\beta)p_j^2
-\{\bar{b}_j^2+\beta(1-\beta)p_j^2\}\right]\\
&=-\frac{\beta(1-\beta)\bar{\lambda}_j^2}{K_j^2}[\bar{b}_j-\beta p_j]^2
=\frac{\alpha}{(1-\alpha)^2}\frac{(p_j+q_j)^2}{K_j^2}\bar{\lambda}_j^2.
\end{split}
\]
This and $p_j^2\bar{R}_j=\bar{b}_j+K_j$ imply
\[
\bar{g}_j=
 \frac{\alpha}{(1-\alpha)^2}\frac{(p_j+q_j)^2}{K_j^2}\bar{\lambda}_j^2
 -(\bar{b}_j+K_j).
\]
Since $(1-\alpha)(1-\beta)=1$, it follows that
\[
(1-\alpha)\bar{b}_j=(1-\alpha)[(\beta-1)p_j-q_j]=-p_j+(\alpha-1)q_j
=q_j\alpha-(p_j+q_j).
\]
Also,
\[
K_j^2=(p_j+q_j)^2-\beta p_j(p_j+2q_j)
=(p_j+q_j)^2+\frac{\alpha}{1-\alpha} p_j(p_j+2q_j).
\]
Combining, we obtain \eqref{eq:3.13}.
\end{proof}

Recall $\xi(t)$ from \eqref{eq:2.3}. 
Taking into account (\ref{eq:2.21}), 
we consider 
$\hat{\pi}=(\hat{\pi}(t))_{t\ge 0}\in \mathcal{A}$ defined by
\begin{equation}
\hat{\pi}(t):=(\sigma')^{-1}(t)
\left[(1-\beta)\{\lambda(t)-\xi(t)\}
-p\bar{R}\xi(t)+p\bar{v}\right],
\quad t\ge 0,
\label{eq:3.14}
\end{equation}
where 
$p:=\mathrm{diag}(p_1,\dots,p_n)$ as in Section \ref{sec:2}, 
and $\bar{R}:=\mathrm{diag}(\bar{R}_1,\dots,\bar{R}_n)$, 
$\bar{v}:=(\bar{v}_1,\dots,\bar{v}_n)'$.

Recall that we have assumed \eqref{eq:2.1}, \eqref{eq:3.1} 
and \eqref{eq:3.2}. Recall also $\alpha^*$ from \eqref{eq:3.3} with 
\eqref{eq:3.4}. 
Here is the solution to Problem P2.

\begin{thm}
\label{thm:3.4}
Let $\alpha^*<\alpha<1$, $\alpha\ne 0$. 
Then $\hat{\pi}$ is an optimal strategy for Problem P2 
with limit in {\rm \eqref{eq:3.6}}. 
The optimal growth rate $J(\alpha)$ in {\rm (P2)} is given by
\begin{equation}
J(\alpha)=\bar{r}+\frac{1}{2\alpha}\sum_{j=1}^nF_j(\alpha)
+\frac{1}{2\alpha}\sum_{j=1}^{n}G_j(\alpha),
\label{eq:3.15}
\end{equation}
where $F_j$'s and $G_j$'s are as in {\rm \eqref{eq:3.10}} and 
{\rm \eqref{eq:3.11}}, respectively.
\end{thm}

\begin{proof}
For simplicity, we put $X(t):=X^{x,\hat{\pi}}(t)$. 
For $\tilde{J}(\alpha)$ in {\rm \eqref{eq:3.5}} and $J^*(\alpha)$ 
in {\rm \eqref{eq:3.6}}, 
we claim $J^*(\alpha)=\tilde{J}(\alpha)$, that is,
\begin{equation}
 \lim_{T\to\infty}\frac{1}{\alpha T}\log E[X^{\alpha}(T)]=\tilde{J}(\alpha).
\label{eq:3.16}
\end{equation}
As mentioned before, \eqref{eq:3.7} and \eqref{eq:3.16} imply
that $\hat{\pi}$ is an optimizer 
for Problem P2. The equality \eqref{eq:3.15} follows from 
this and \eqref{eq:3.13}

We complete the proof of the theorem by proving \eqref{eq:3.16}.

\noindent {\it Step}\/ 1.\ We calculate $E[X^{\alpha}(T)]$. 
Define the $\mathbf{R}$-valued martingale $L(t)$ by
\[
L(t):=\exp\left[\alpha\int_0^t\{\sigma'(s)\hat{\pi}(s)\}'dB(s)
-\frac{\alpha^2}{2}\int_0^t\Vert \sigma'(s)\hat{\pi}(s)
 \Vert^2ds\right], \quad t\ge 0.
\]
From (\ref{eq:2.4}), we have 
$X^{\alpha}(t)=[xS_0(t)]^{\alpha}L(t)\exp[\int_{0}^{t}N(s)ds]$ 
for $t\ge 0$, 
where
\[
\begin{split}
&N(t):=\alpha\{\sigma'(t)\hat{\pi}(t)\}'\left[\lambda(t)-\xi(t)+
 \frac{1}{2}(\alpha-1)\sigma'(t)\hat{\pi}(t)\right]\\
&=\frac{\alpha(1-\alpha)}{2}\{\sigma'(t)\hat{\pi}(t)\}'
\left[
(1-\beta)\{\lambda(t)-\xi(t)\}+p\bar{R}\xi(t)-p\bar{v}
\right]\\
&=-\frac{\beta}{2}
\left[\lambda(t)-\xi(t)-(1-\alpha)\{p\bar{R}\xi(t)-p\bar{v}\}\right]'\\
&\qquad\qquad\qquad\qquad 
\cdot \left[\lambda(t)-\xi(t)+(1-\alpha)\{p\bar{R}\xi(t)-p\bar{v}\}\right]\\
&=-\frac{\beta}{2}
\left[
\{\lambda(t)-\xi(t)\}'\{\lambda(t)-\xi(t)\}
-(1-\alpha)^2\{p\bar{R}\xi(t)-p\bar{v}\}'\{p\bar{R}\xi(t)-p\bar{v}\}
\right].
\end{split}
\]
Notice that we have used 
$(1-\alpha)(1-\beta)=1$, $\alpha(1-\beta)=-\beta$. 
We write
\[
N(t)= -\frac{1}{2}\xi'(t)Q\xi(t)-h'(t)\xi(t)+\frac{1}{2}
\sum\nolimits_{j=1}^nu_j(t),
\]
where
\[
u_j(t):=\alpha(\alpha-1)\left[p_j^2\bar{v}_j^2
    -(1-\beta)^2\lambda_j^2(t)\right],\quad t\ge 0,\quad j=1,\dots,n,
\]
and $Q=\mathrm{diag}(Q_1,\dots,Q_n)$, 
$h(t)=h(t;T)=(h_1(t;T),\dots,h_j(t;T))'$ with
\begin{align*}
&Q_j:=\beta\left[1-(1-\alpha)^2p_j^2\bar{R}_j^2\right],\quad 
t\ge 0,\quad j=1,\dots,n,\\
&h_j(t):=\alpha(\alpha-1)p_j^2\bar{R}_j\bar{v}_j-\beta\lambda_j(t),
\quad t\ge 0,\quad j=1,\dots,n.
\end{align*}
Therefore,
\begin{equation}
\begin{split}
 E[X^{\alpha}(T)]
  &=[xS_0(T)]^{\alpha}
  \exp\left[\frac{1}{2}\sum\nolimits_{j=1}^n \int_0^T u_j(t)dt\right]\\
 &\quad\times\bar{E}\left[\exp\left\{-\int_0^T
 \left(\frac{1}{2}\xi'(t)Q\xi(t)+h'(t)\xi(t)\right)dt
  \right\}\right],
\end{split}
\label{eq:3.17}
\end{equation}
where $\bar{E}$ denotes the expectation with respect to 
the probability measure $\bar{P}$ on $(\Omega,\mathcal{F}_T)$ such that 
$d\bar{P}/dP=L(T)$.

\noindent {\it Step}\/ 2.\ We continue the calculation of $E[X^{\alpha}(T)]$. 
We are about to apply Theorem \ref{thm:A1} in Appendix A to \eqref{eq:3.17}. 
Write $\bar{B}(t):=B(t)-\alpha\int_0^t\sigma'(s)\hat{\pi}(s)ds$ for $t\ge 0$. 
Then $\bar{B}(t)$ is an $\mathbf{R}^n$-valued standard Brownian motion under 
$\bar{P}$. By \eqref{eq:2.6}, the process $\xi(t)$ evolves according to 
the $n$-dimensional stochastic differential equation
\begin{equation}
 d\xi(t)=[\gamma(t)+d(t)\xi(t)]dt + l(t)d\bar{B}(t),\quad t\ge 0, 
\label{eq:3.18}
\end{equation}
where 
$d(t)=\mathrm{diag}(d_1(t),\dots,d_n(t))$, 
$\gamma(t)=\mathrm{diag}(\gamma_1(t),\dots,\gamma_n(t))$ with
\begin{align*}
&d_j(t):=b_j(t)-\alpha p_j\bar{R}_jl_j(t),\quad 
t\ge 0,\quad j=1,\dots,n,\\
&\gamma_j(t):=\rho_j(t)+\alpha p_jl_j(t)\bar{v}_j,
\quad t\ge 0,\quad j=1,\dots,n.
\end{align*}
For $j=1,\dots,n$, let $U_j(t)\equiv U_j(t;T)$ be the unique solution 
to the one-dimensional backward Riccati equation
\begin{equation}
 \dot{U}_j(t)-l_j^2(t)U_j^2(t)+2d_j(t)U_j(t)+Q_j=0,
 \quad 0\le t\le T,\quad U_j(T)=0
 \label{eq:3.19}
\end{equation}
in the sense of Lemma \ref{lem:3.5} below, and 
let $m_j(t)\equiv m_j(t;T)$ be the solution to the 
one-dimensional linear equation
\begin{equation}
\begin{split}
 &\dot{m}_j(t)+[d_j(t)-l_j^2(t)U_j(t;T)]m_j(t)
 -h_j(t)-U_j(t;T)\gamma_j(t)=0\\
&\quad 0\le t\le T, \quad m_j(T)=0.
 \end{split}
\label{eq:3.20}
\end{equation}
Then, from \eqref{eq:3.17}--\eqref{eq:3.20} and Theorem \ref{thm:A1}, we obtain
\begin{equation}
E[X^{\alpha}(T)]
  =[xS_0(T)]^{\alpha}
  \exp\left[\frac{1}{2}\sum\nolimits_{j=1}^{n}\int_0^T f_j(t;T)dt\right],
\label{eq:3.21}
\end{equation}
where, for $(t,T)\in\Delta$ and $j=1,\dots,n$,
\[
f_j(t;T):=l_j^2(t)m_j^2(t;T)+2\gamma_j(t)m_j(t;T)
-l_j^2(t)U_j(t;T)+u_j(t).
\]

\noindent {\it Step}\/ 3.\ We compute the limit 
$J^*(\alpha)$ in \eqref{eq:3.6}. Let $j\in\{1,\dots,n\}$. 
Write
\[
\bar{d}_j:=\bar{b}_j-\alpha p_j^2\bar{R}_j.
\]
Then $d_j(t)$ converges to $\bar{d}_j$, as $t\to\infty$, exponentially fast. 
Now
\begin{align*}
\bar{d}_j^2 + p_j^2Q_j&= (\bar{b}_j-\alpha p_j^2 \bar{R}_j)^2 
 + p_j^2\beta\left[1-(1-\alpha)^2p_j^2\bar{R}_j^2\right] \\
&=\bar{b}_j^2 -2\alpha\bar{b}_j(\bar{b}_j+K_j)+\alpha^2(\bar{b}_j+K_j)^2
 + p_j^2\beta-\alpha(\alpha-1)(\bar{b}_j+K_j)^2 \\
&=(1-\alpha)\bar{b}_j^2 +\alpha K_j^2 +p_j^2\beta
=\bar{b}_j^2 +p_j^2\beta(1-\beta),
\end{align*} 
which implies
\begin{equation}
\sqrt{\bar{d}_j^2 + p_j^2Q_j}=K_j>0.
\label{eq:3.22}
\end{equation}
Thus we may write $\bar{U}_j$ for the larger (resp.\ unique) solution of 
the following equation when $p_j>0$ (resp.\ $p_j=0$):
\begin{equation}
p_j^2x^2 -2\bar{d}_jx-Q_j=0.
\label{eq:3.23}
\end{equation}
From (\ref{eq:3.22}), we also see that $\bar{d}_j-p_j^2\bar{U}_j=-K_j$. 
Let $\bar{m}_j$ be the solution to
\begin{equation}
\left(\bar{d}_j-p_j^2\bar{U}_j\right)\bar{m}_j-\bar{h}_j
-\bar{U}_j\bar{\gamma}_j=0,
\label{eq:3.24}
\end{equation}
where
\[
\bar{h}_j:=\alpha(\alpha-1)p_j^2\bar{R}_j\bar{v}_j-\beta\bar{\lambda}_j,
\quad 
\bar{\gamma}_j:=\bar{\rho}_j+\alpha p_j^2\bar{v}_j.
\]

By \eqref{eq:3.21}, we have
\[
 \frac{1}{\alpha T}\log E[X^{\alpha}(T)] 
 = \frac{\log x}{T}+\frac{1}{T}\int_0^Tr(t)dt
 +\frac{1}{2\alpha}\sum_{j=1}^n\frac{1}{T}\int_0^T f_j(t;T)dt.
\]
However, Propositions \ref{prop:3.6} and \ref{prop:3.7} below 
imply
\[
\lim_{T\to\infty}\frac{1}{T}\int_{0}^{T}f_j(t;T)dt=\bar{f}_j
\]
with
\[
\bar{f}_j:=p_j^2\bar{m}_j^2+2\bar{\gamma}_j\bar{m_j}
-p_j^2\bar{U}_j+\alpha(\alpha-1)\left[p_j^2\bar{v}_j^2
    -(1-\beta)^2\bar{\lambda}_j^2\right],
\]
so that
\begin{equation}
J^*(\alpha)=\bar{r}+\frac{1}{2\alpha}\sum_{j=1}^{n}\bar{f}_j.
\label{eq:3.25}
\end{equation}

\noindent {\it Step}\/ 4.\ 
Here we show that in fact \eqref{eq:3.16} holds. 
First,
\[
p_j^2\bar{U}_j=\bar{d}_j+K_j=\bar{b}_j-\alpha p_j^2\bar{R}_j+K_j
=p_j^2\bar{R}_j-\alpha p_j^2\bar{R}_j=(1-\alpha)p_j^2\bar{R}_j,
\]
whence $\bar{U}_j=(1-\alpha)\bar{R}_j$ 
(which we can directly check when $p_j=0$). 
Next, 
\[
\begin{aligned}
\bar{h}_j+\bar{U}_j\bar{\gamma}_j
&=\alpha(\alpha-1)p_j^2\bar{R}_j\bar{v}_j-\beta\bar{\lambda}_j
+(1-\alpha)\bar{R}_j[-\beta p_j\bar{\lambda}_j+\alpha p_j^2\bar{v}_j]\\
&=\bar{\lambda}_j(-\beta +\alpha p_j\bar{R}_j)
=-(1-\alpha)\beta \bar{\lambda}_j[(1-\beta)+p_j\bar{R}_j],
\end{aligned}
\]
so that
\[
\bar{m}_j=\frac{(1-\alpha)}{K_j}\beta \bar{\lambda}_j[(1-\beta)+p_j\bar{R}_j]
=(1-\alpha)\bar{v}_j.
\]
Therefore, 
\[
\begin{split}
\bar{f}_j
&=(1-\alpha)^2p_j^2\bar{v}^2_j
+2(1-\alpha)(\bar{\rho}_j+\alpha p_j^2\bar{v}_j)\bar{v}_j
-(1-\alpha)p_j^2\bar{R}_j\\
&\qquad\qquad\qquad\qquad\qquad\qquad\qquad\qquad 
+\alpha(\alpha-1)\left[p_j^2\bar{v}_j^2
    -(1-\beta)^2\bar{\lambda}_j^2\right]\\
&=(1-\alpha)
\left[
p_j^2\bar{v}^2_j+2\bar{\rho}_j\bar{v}_j-p_j^2\bar{R}_j
 -\beta(1-\beta)\bar{\lambda}_j^2
\right]=(1-\alpha)\bar{g}_j.
\end{split}
\]
From \eqref{eq:3.12}, \eqref{eq:3.25} and this, we obtain 
$J^*(\alpha)=\tilde{J}(\alpha)$ or (\ref{eq:3.16}), as desired.
\end{proof}

In the proof above, we needed the following results.

\begin{lem}\label{lem:3.5}Let $j\in\{1,\dots,n\}$.
\begin{enumerate}
\item If $p_j=0$, then {\rm \eqref{eq:3.19}} has a unique 
solution $U_j(t)\equiv U_j(t;T)$.
\item If $p_j>0$ and $\alpha_j^*<\alpha<0$, then {\rm \eqref{eq:3.19}} has 
a unique nonnegative solution $U_j(t)\equiv U_j(t;T)$.
\item If $p_j>0$ and $0<\alpha<1$, then 
{\rm \eqref{eq:3.19}} has a unique solution $U_j(t)\equiv U_j(t;T)$ such that 
$U_j(t;T)\ge (1-\alpha)R_j(t;T)$ for $t\in [0,T]$, where 
$R_j(t)\equiv R_j(t;T)$ is the solution to {\rm \eqref{eq:2.10}} in the sense 
of Lemma {\rm \ref{lem:2.1} (iii)}.
\end{enumerate}
\end{lem}

\begin{proof}(i)\ When $p_j=0$, \eqref{eq:3.19} is linear, whence it 
has a unique solution.

(ii)\ For $p_j>0$ and $\alpha<0$, we put 
$f(x)=p_j^2x^2-2\bar{b}_jx-\beta(1-\beta)$. 
Since $\bar{b}_j<0$ and $\beta(1-\beta)>0$, 
the larger solution $\bar{R}_j$ to $f(x)=0$ satisfies 
$p_j^2\bar{R}_j^2 <(1-\beta)^2$ if and only if 
$f((1-\beta)/p_j)>0$. However, 
this is equivalent to $-3p_j-2q_j<(p_j-2q_j)\alpha$. Thus, if $p_j>0$ 
and $\alpha_j^*<\alpha<0$, then $p^2_j\bar{R}^2_j <(1-\beta)^2$ or 
$Q_j>0$, so that the Riccati equation 
\eqref{eq:3.19} has a unique nonnegative solution.

(iii)\ Suppose $p_j>0$ and $0<\alpha<1$. 
For the solution $R_j(t)\equiv R_j(t;T)$ to \eqref{eq:2.10} in the sense of 
Lemma \ref{lem:2.1} (iii), we consider
\[
P_j(t):=\frac{U_j(t)}{1-\alpha}-R_j(t). 
\]
Let $d_j(t)$ be as above. Then, \eqref{eq:3.19} becomes
\[
\begin{split}
&\dot{P}_j(t)-(1-\alpha)l_j^2(t)P^2_j(t) 
 -2\left[(1-\alpha)l_j^2(t)R_j(t)-d_j(t)\right]P_j(t)\\
&\qquad\qquad\qquad\qquad\qquad\qquad\qquad
+\alpha[l_j(t)R_j(t)-p_j\bar{R}_j]^2=0,\quad 
0\le t\le T,
\end{split}
\]
with $P_j(T)=0$. Since $(1-\alpha)l_j^2(t)>0$ and 
$\alpha[l_j(t)R_j(t)-p_j\bar{R}_j]^2>0$, 
this Riccati equation has a unique nonnegative solution. 
Thus the assertion follows.
\end{proof}

\begin{prop}\label{prop:3.6}Let $\alpha^*<\alpha<1$, $\alpha\ne 0$, and 
$j\in\{1,\dots,n\}$. Let $U_j(t;T)$ be the unique solution to 
{\rm \eqref{eq:3.19}} in the sense of Lemma {\rm \ref{lem:3.5}}, 
and let $\bar{U}_j$ be the larger (resp.\ unique) solution to 
{\rm \eqref{eq:3.23}} when $p_j>0$ (resp.\ $p_j=0$). 
Then
\begin{enumerate}
\item $U_j(t;T)$ is bounded in $\Delta$.
\item $\lim_{T-t\to\infty,\ t\to\infty}U_j(t;T)=\bar{U}_j$.
\item For $\delta, \epsilon\in (0,\infty)$ such that 
$\delta+\epsilon<1$,
\[
\lim_{T\to\infty}\sup_{\delta T\le t\le (1-\epsilon)T}
 |U_j(t;T)-\bar{U}_j|=0.
\]
\end{enumerate}
\end{prop}

\begin{proof}
We assume $0<\alpha<1$ and $p_j>0$. Since 
$U_j(t;T)\ge (1-\alpha)R_j(t;T)$ in $\Delta$ and 
$R_j(t;T)$ is bounded from below 
by Proposition \ref{prop:3.1}, so is $U_j(t;T)$. 
Let $N_j(t)\equiv N_j(t;T)$ be the solution to the linear equation
\[
 \dot{N}_j(t)+2[d_j(t)-l_j^2(t)\bar{U}_j]N_j(t)+Q_j+l_j^2(t)\bar{U}_j^2=0, 
\quad 0\le t\le T,\quad N_j(T)=0.
\]
By \eqref{eq:3.22}, 
$d_j(t)-l_j^2(t)\bar{U}_j\to \bar{d}_j-p_j^2\bar{U}_j=-K_j<0$ 
as $t\to\infty$, 
so that $N_j(t;T)$ is bounded from above in $\Delta$. 
Since $N_j(T)-U_j(T)=0$ and
\[
 [\dot{N}_j(t)-\dot{U}_j(t)]+2[d_j(t)-l_j^2(t)\bar{U}_j][N_j(t)-U_j(t)]
 =-l_j^2(t)\left[U_j(t)-\bar U_j\right]^2\le 0,
\]
we have, as in the proof of Proposition \ref{prop:3.1}, 
$U_j(t;T)\le N_j(t;T)$ in $\Delta$. Thus $U_j(t;T)$ is also bounded from 
above in $\Delta$. Combining, $U_j(t;T)$ is bounded in $\Delta$. 
The rest of the proof is similar to that of Proposition 
\ref{prop:3.1}, whence we omit it.
\end{proof}

\begin{prop}\label{prop:3.7}
Let $\alpha^*<\alpha<1$, $\alpha\ne 0$, and 
$j\in\{1,\dots,n\}$. Let $m_j(t;T)$ and $\bar{m}_j$ be the solutions to 
{\rm \eqref{eq:3.20}} and {\rm \eqref{eq:3.24}}, respectively. Then
\begin{enumerate}
\item $m_j(t;T)$ is bounded in $\Delta$.
\item $\lim_{T-t\to\infty,\ t\to\infty}m_j(t;T)=\bar{m}_j$.
\item For $\delta, \epsilon\in (0,\infty)$ 
such that $\delta+\epsilon<1$, 
\[
\lim_{T\to\infty}\sup_{\delta T\le t\le (1-\epsilon)T}
\vert m_j(t;T)-\bar{m}_j\vert=0.
\]
\end{enumerate}
\end{prop}

The proof of Proposition \ref{prop:3.7} 
is similar to that of Proposition \ref{prop:3.2}; so 
we omit it.

\begin{rem}
\label{rem:3.8}
We note that the proof of Lemma \ref{lem:3.5} (iii) is still valid 
under \eqref{eq:1.4} if there were a solution $R_j(t)\equiv R_j(t;T)$ to 
\eqref{eq:2.10}. This implies that, to prove an analogue of 
Theorem \ref{thm:3.4} with $0<\alpha<1$, which is relevant to Problem P3, 
for \eqref{eq:1.4}, one may show 
the existence of such $R_j(t)$ when $0<\alpha<1$. We did not succeed in 
such an extension to Lemma \ref{lem:2.1} (iii) (see Remark \ref{eq:2.6}).
\end{rem}


\section{Large deviations probability control}\label{sec:4}

In this section, we study the large deviations probability control problem 
P3 for the market model $\mathcal{M}$. 
Throughout this section, we assume (\ref{eq:2.1}), (\ref{eq:3.1}), 
(\ref{eq:3.2}) and
\begin{equation}
\mbox{either $\bar{\lambda}\ne (0,\dots,0)'$ or 
$(p_1,\dots,p_n)\ne (0,\dots,0)$.}
\label{eq:4.1}
\end{equation}

For $x\in (0,\infty)$ and $\pi\in\mathcal{A}$, 
let $L^{x,\pi}(T)$ be the growth rate defined by
\[
L^{x,\pi}(T):=\frac{\log X^{x,\pi}(T)}{T},\quad T>0.
\]
We have $P\left(L^{x,\pi}(T)\ge c\right)
=P\left(X^{x,\pi}(T)\ge e^{cT}\right)$. 
Following Pham \cite{P1,P2}, 
we consider the optimal logarithmic moment generating function
\[
\Lambda(\alpha):=\sup_{\pi\in\mathcal{A}} 
\limsup_{T\to\infty} \log E[\exp(\alpha TL^{x,\pi}(T))],\quad 0<\alpha<1.
\]
Since $\Lambda(\alpha)=\alpha J(\alpha)$ for $\alpha\in (0,1)$, 
it follows from Theorem \ref{thm:3.4} that
\[
\Lambda(\alpha)=\bar{r} \alpha+\frac{1}{2}\sum_{j=1}^{n}F_{j}(\alpha)
+\frac{1}{2}\sum_{j=1}^{n}G_{j}(\alpha),\quad 0<\alpha<1,
\]
where $F_j$'s and $G_j$'s are as in {\rm \eqref{eq:3.10}} and 
{\rm \eqref{eq:3.11}}, respectively.

\begin{prop}
\label{prop:4.1}
We have $(d\Lambda/d\alpha)(0+)=\bar{c}$ and 
$\lim_{\alpha\uparrow 1}(d\Lambda/d\alpha)(\alpha)=\infty$, where
\[
\bar{c}
:=\bar{r}+\frac{1}{4}\sum_{j=1}^n\frac{p_j^2}{p_j+q_j}
+\frac{1}{2}\Vert \bar{\lambda}\Vert^2.
\]
\end{prop}

\begin{proof}
For $0<\alpha<1$, $\dot{F}_j(\alpha)$ is equal to
\[
\frac{(p_j+q_j)^2\bar{\lambda}_j^2}{\left[(1-\alpha)(p_j+q_j)^2
+\alpha p_j(p_j+2q_j)\right]}
+\frac{(p_j+q_j)^2\bar{\lambda}_j^2q_j^2\alpha}
{\left[(1-\alpha)(p_j+q_j)^2+\alpha p_j(p_j+2q_j)\right]^2}.
\]
From this, $\dot{F}_j(0+)=\bar{\lambda}_j^2$. 
This also shows that
\[
\frac{dF_j}{d\alpha}(\alpha)\sim \bar{\lambda}_j^2(1-\alpha)^{-2},\quad 
\alpha\uparrow 1
\]
if $p_j=0$ and $\bar{\lambda}_j\ne 0$. 
On the other hand, for $0<\alpha<1$,
\[
\begin{split}
&\frac{dG_j}{d\alpha}(\alpha)=
-q_j+
\frac{(1-\alpha)^{-1/2}}{2}
\left[(1-\alpha)(p_j+q_j)^2+\alpha p_j(p_j+2q_j)\right]^{1/2}\\
&\qquad\qquad\qquad\qquad 
+
\frac{q_j^2(1-\alpha)^{1/2}}{2\left[(1-\alpha)(p_j+q_j)^2
+\alpha p_j(p_j+2q_j)\right]^{1/2}}.
\end{split}
\]
This gives $(dG_j/d\alpha)(0+)=p_j^2/[2(p_j+q_j)]$. This 
also yields
\[
\frac{dG_j}{d\alpha}(\alpha)\sim \frac{\sqrt{p_j(p_j+2q_j)}}{2}
(1-\alpha)^{-1/2},\quad 
\alpha\uparrow 1
\]
if $p_j>0$. Thus the proposition follows.
\end{proof}

\begin{rem}
\label{rem:4.2}
From the proof of Proposition \ref{prop:4.1}, we see that
\[
\frac{d\Lambda}{d\alpha}(\alpha)\sim 
\frac{(1-\alpha)^{-1/2}}{4}
\sum\nolimits_{j=1}^n\sqrt{p_j(p_j+2q_j)}
,\quad \alpha\uparrow 1
\]
if $p_j>0$ for all $j=1,\dots,n$, otherwise
\[
\frac{d\Lambda}{d\alpha}(\alpha)\sim 
\frac{(1-\alpha)^{-2}}{2}\sum\nolimits_{1\le j\le n \atop p_j=0}
\bar{\lambda}_j^2,\quad \alpha\uparrow 1.
\]
\end{rem}

For $\alpha\in (0,1)$, we denote by $\hat{\pi}(t;\alpha)$ the optimal 
strategy $\hat{\pi}(t)$ in (\ref{eq:3.14}). 
Recall $I(c)$ from (P3). 
From Theorem \ref{thm:3.4}, Proposition \ref{prop:4.1}, 
and Pham \cite[Theorem 3.1]{P1}, we immediately obtain 
the following solution to Problem P3:

\begin{thm}
\label{thm:4.3}
We have
\[
I(c)=-\sup_{\alpha\in (0,1)}\left[\alpha c-\Lambda(\alpha)\right],
\quad c\in\mathbf{R}.
\]
Moreover, if $\alpha(d)\in (0,1)$ is such that 
$\dot{\Lambda}(\alpha(d))=d\in (\bar{c},\infty)$, then, 
for $c\ge \bar{c}$, the sequence of strategies
\[
\hat{\pi}^m(t):=\hat{\pi}(t;\alpha(c+\tfrac{1}{m}))
\]
is nearly optimal in the sense that
\[
\lim_{m\to\infty} \limsup_{T\to\infty}\frac{1}{T}\log 
P\left(X^{x,\hat{\pi}^m}(T)\ge e^{cT}\right)=I(c),\quad c\ge\bar{c}.
\]
\end{thm}

\begin{rem}\label{rem:4.4}
Theorem 3.1 in Pham \cite{P1} is stated for a 
model different from $\mathcal{M}$ but the arguments there are 
so general that we can prove Theorem \ref{thm:4.3} in the same way.
\end{rem}

We turn to the problem of deriving an optimal strategy, 
rather than a nearly optimal sequence, for the problem (P3)
when $c<\bar{c}$. 
We define $\pi_0\in\mathcal{A}$ by
\[
\hat{\pi}_{0}(t):= (\sigma')^{-1}(t)
\left[\lambda(t)-\xi(t)\right], \quad t\ge 0,
\]
where recall $\xi(t)$ from \eqref{eq:2.3}. 
From (\ref{eq:2.4}), 
\[
\begin{split}
L^{x,\hat{\pi}_0}(T)
&=\frac{\log x}{T} +\frac{1}{T}\int_0^T r(t)dt
 +\frac{1}{2T}\int_0^T\left\Vert\lambda(t)-\xi(t)\right\Vert^2dt\\
&\qquad\qquad\qquad\qquad\qquad\qquad\qquad 
+\frac{1}{T}\int_0^T\left[\lambda(t) - \xi(t)\right]'dB(t).
\end{split}
\]

\begin{prop}
\label{prop:4.5}
The rate $L^{x,\hat{\pi}_0}(T)$ converges to $\bar{c}$, 
as $T\to\infty$, in probability.
\end{prop}

\begin{proof}
In this proof, we denote by $C$ positive constants, which may not be 
necessarily equal.

For $j=1,\dots,n$, we write
\[
\lambda_j(t)-\xi_j(t)
=[\lambda_j(t)-\bar{\lambda}_j]+[\bar{\lambda}_j-p_jK(t)]+N(t),
\]
where $K(t)=\int_0^te^{-(p_j+q_j)(t-s)}dB_j(s)$ and 
$N(t)=\int_0^te^{-(p_j+q_j)(t-s)}f(s)dB_j(s)$
with
\[
f(s)=\frac{2p_j^2q_j}{(2q_j+p_j)^2e^{2q_js}-p_j^2}.  
\]
The process $K(t)$, the dynamics of which are given by 
\[
dK(t)=-(p_j+q_j)K(t)dt+dB_j(t),
\]
is a positively recurrent one-dimensional diffusion process with 
speed measure $m(dx)=2e^{-(p_j+q_j)x^2}dx$. By the ergodic theorem 
(cf.\ Rogers and Williams \cite[v.53]{RW}), we have
\[
\lim_{T\to\infty}\frac{1}{T}\int_0^T[\bar{\lambda}_j-p_jK(t)]^2dt
=\int_{-\infty}^{\infty}(\bar{\lambda}_j-p_jy)^2\nu(dy)
=\bar{\lambda}_j^2+\frac{p_j^2}{2(p_j+q_j)}\quad \mbox{a.s.,}
\]
where $\nu(dy)$ is the 
Gaussian measure with mean $0$ and variance $1/[2(p_j+q_j)]$.

Since $0\le f(s)\le Ce^{-2q_js}$, we have
\[
E\left[N^2(t)\right]\le C\int_0^te^{-2q_j(t+s)}ds
\le Ce^{-2q_jt}, \quad t\ge 0.
\]
Also, $E[K^2(t)]\le C$ for $t\ge 0$. Therefore, 
\[
\frac{1}{T}\int_0^TE\left[\vert\{\bar{\lambda}_j-p_jK(t)\}N(t)\vert\right]
dt
\le \frac{C}{T}\int_0^TE[N^2(t)]^{1/2}dt \to 0,\quad T\to \infty.
\]
Similarly, 
\[
\begin{split}
&\lim_{T\to\infty}\frac{1}{T}\int_0^T[\lambda_j(t)-\bar{\lambda}_j]^2dt
=\lim_{T\to\infty}\frac{1}{T}\int_0^TE\left[(\lambda_j(t)-\bar{\lambda}_j)
(\bar{\lambda}_j-p_jK(t))\right]dt\\
&\quad=\lim_{T\to\infty}\frac{1}{T}\int_0^TE\left[N^2(t)\right]dt
=\lim_{T\to\infty}\frac{1}{T}\int_0^TE\left[(\lambda_j(t)-\bar{\lambda}_j)N(t)
\right]dt
=0.
\end{split}
\]
Combining, 
\[
\frac{1}{T}\int_0^T\left[\lambda_j(t)-\xi_j(t)\right]^2dt\to
\bar{\lambda}_j^2+\frac{p_j^2}{2(p_j+q_j)}, \quad T\to\infty,\quad 
\mbox{\rm in probability.} 
\]

Finally, for $j=1,\dots,n$ and $t\ge 0$, 
\[
E\left[\left\{\lambda_j(t)-\xi_j(t)\right\}^2\right]
\le 2\lambda^2_j(t)+2E\left[\xi_j^2(t)\right]
\le C\left[1+\int_0^tl^2_j(t,s)ds\right]\le C,
\]
so that $(1/T)\int_0^T[\lambda_j(t)-\xi_j(t)]dB_j(t)\to 0$, 
as $T\to\infty$, in $L^2(\Omega)$, 
whence in probability. 
Thus the proposition follows.
\end{proof}

\begin{thm}\label{thm:4.6}
For $c<\bar{c}$, $\hat{\pi}_0$ is optimal for Problem P3 
with limit
\[
\lim_{T\to\infty}\frac{1}{T}\log 
P\left(X^{x,\hat{\pi}_0}(T)\ge e^{cT}\right)=I(c),\quad c<\bar{c}.
\]
\end{thm}

\begin{proof}
Proposition \ref{prop:4.5} implies 
$\lim_{T\to\infty}P\left(L^{x,\hat{\pi}_0}(T)\ge c\right)=1$ 
for $c<\bar{c}$, so that
\[
\lim_{T\to\infty}\frac{1}{T}\log 
 P\left(L^{x,\hat{\pi}_0}(T)\ge c\right)=0\ge 
 \sup_{\pi\in\mathcal{A}}\lim_{T\to\infty}\frac{1}{T}\log 
 P\left(L^{x,\pi}(T)\ge c\right),\quad 
 c<\bar{c}.  
\]
Thus $\hat{\pi}_0$ is optimal if $c<\bar{c}$.
\end{proof}

\begin{rem}\label{rem:4.7}
From Theorem 10.1 in Karatzas and Shreve \cite[Chapter 3]{KS}, 
we see that 
$\hat{\pi}_0$ is the log-optimal or growth optimal strategy in 
the sense that
\[
\sup_{\pi\in\mathcal{A}}\limsup_{T\to\infty}\frac{1}{T}
 \log X^{x,\pi}(T) = \limsup_{T\to\infty}\frac{1}{T}
 \log X^{x,\hat{\pi}_{0}}(T) \quad\mbox{a.s.} 
\]
We also find that $\lim_{\alpha\downarrow 0}\hat{\pi}(t;\alpha)
=\hat{\pi}_{0}(t)$ a.s.\ for 
$t\ge 0$. 
\end{rem}


\appendix

\section{A Cameron--Martin type formula}

In this appendix, we prove a generalization of the Cameron--Martin formula that we need in the proofs of Proposition \ref{prop:2.2} and 
Theorem \ref{thm:3.4}. We refer to Myers \cite{M} for earlier work.

Let $T\in (0,\infty)$  and let $\mathbf{R}^{n\times n}$ be the set of 
$n\times n$ real matrices. 
We say that $A:[0,T]\to\mathbf{R}^{n\times n}$ is symmetric if $A(t)$ is a 
symmetric matrix for all $t\in [0,T]$. 
Let $(\Omega,\mathcal{F},P)$ be the underlying 
complete probability space equipped with filtration 
$(\mathcal{F}_t)_{0\le t\le T}$ satisfying 
the usual conditions. 
We assume that the $\mathbf{R}^n$-valued process 
$\xi(t)$ satisfies the $n$-dimensional stochastic differential equation
\[
 d\xi(t)=[a(t)+b(t)\xi(t)]dt + c(t)dB(t),\quad 0\le t\le T,
\]
where $B(t)$ is an $\mathbf{R}^n$-valued standard $(\mathcal{F}_t)$-Brownian 
motion and all the coefficients $a:[0,T]\to\mathbf{R}^n$ and 
$b,c:[0,T]\to\mathbf{R}^{n\times n}$ 
are deterministic, bounded measurable functions.

\begin{thm}
\label{thm:A1}
Let $Q:[0,T]\to\mathbf{R}^{n\times n}$ and $h:[0,T]\to\mathbf{R}^{n}$ be 
deterministic, bounded measurable functions. We assume that $Q$ 
is symmetric. 
We also assume that there exists a bounded symmetric function 
$R:[0,T]\to\mathbf{R}^{n\times n}$ satisfying 
the backward matrix Riccati equation
\begin{equation}
\begin{split}
&\dot{R}(t)-R(t)c(t)c'(t)R(t)+b'(t)R(t)+R(t)b(t)+ Q(t)=0, 
\quad 0\le t\le T,\\
&R(T)=0.
\end{split}
\label{eq:A.1}
\end{equation}
Let $v:[0,T]\to\mathbf{R}^n$ be the solution to the linear equation
\begin{equation}
\begin{split}
&\dot{v}(t)+[b(t)-c(t)c'(t)R(t)]'v(t)-h(t)-R(t)a(t)=0, 
\quad 0\le t\le T,\\
&v(T)=0.
\end{split}
\label{eq:A.2}
\end{equation}
Then, 
for $t\in [0,T]$,
\begin{align*}
 &E\left[\left.\exp\left\{-\int_t^T
 \left(\frac{1}{2}\xi'(s)Q(s)\xi(s)+h'(s)\xi(s)\right)ds
  \right\}\right|\mathcal{F}_t\right] \\
 &= \exp\left[v'(t)\xi(t)-\frac{1}{2}\xi'(t)R(t)\xi(t)\right]\\
 & \qquad\times\exp\left[\frac{1}{2}\int_t^T\{v'(s)c(s)c'(s)v(s)+2a'(s)v(s)
   -\mathrm{tr}(c(s)c'(s)R(s))\}ds\right]. 
\end{align*}
\end{thm}

\begin{proof}
We put $K(t)=c'(t)[v(t)-R(t)\xi(t)]$ for $t\in [0, T]$. 
Then, by the It\^{o} formula,
\begin{align*}
 &d\left[v'(t)\xi(t)-\frac{1}{2}\xi'(t)R(t)\xi(t) \right] \
 -\left[K'(t)dB(t) - \frac{1}{2}\Vert K(t)\Vert^2dt\right] \\
 &\qquad
 = \left[-\frac{1}{2}\xi'(t)\{\dot{R}(t)+R(t)b(t)+b'(t)R(t)
 -R(t)c(t)c'(t)R(t)\}\xi(t) \right.\\
 &\qquad\qquad\quad\ \ 
 +\{\dot{v}(t)+(b(t)-c(t)c'(t)R(t))'v(t)-R(t)a(t)\}'
 \xi(t) \\
 &\qquad\qquad\qquad\quad\ \ 
 +\left.\frac{1}{2}\{v'(t)c(t)c'(t)v(t)+2a'(t)v(t)
 -\mathrm{tr}(c(t)c'(t)R(t))\}\right]dt\\
 &\qquad =\left[\frac{1}{2}\xi'(t)Q(t)\xi(t)+h'(t)\xi(t)\right.\\
 &\qquad\qquad\qquad\qquad
 \left.+ \frac{1}{2}\{v'(t)c(t)c'(t)v(t)+2a'(t)v(t)
 -\mathrm{tr}(c(t)c'(t)R(t))\}
 \right]dt.
\end{align*}
Therefore, $\int_t^T K'(s)dB_s-\frac{1}{2}\int_t^T \Vert K(s)\Vert^2ds$ 
is equal to
\begin{align*}
 &\frac{1}{2}\xi'(t)R(t)\xi(t) -v'(t)\xi(t)
 -\int_t^T\left(\frac{1}{2}\xi'(s)Q(s)\xi(s)+h'(s)\xi(s)\right)ds\\
 &\qquad\qquad
 -\frac{1}{2}\int_t^T\left\{v'(s)c(s)c'(s)v(s)+2a'(s)v(s)
 -\mathrm{tr}(c(s)c'(s)R(s))\right\}ds . 
\end{align*}
Since $K(t)$ is a continuous Gaussian process, the process
\[
M(t):=\exp\left\{
\int_0^t K'(s)dB(s)-\frac{1}{2}\int_0^t \Vert K(s)\Vert^2ds
\right\},
\quad 0\le t\le T,
\]
is a martingale (cf.\ Example 3(a) in \cite[Section 6.2]{LS}). Thus
\[
 E\left[\left.\exp\left\{\int_t^T K'(s)dB(s)-\frac{1}{2}
  \int_t^T \Vert K(s)\Vert^2ds\right\}\right|\mathcal{F}_t\right]=1.
\]
Combining, we obtain the theorem.
\end{proof}


\section{Asymptotics for a solution to Riccati equation}

Here we summarize the results on the asymptotics for 
a solution to Riccati or linear equation that we need 
in Section \ref{sec:3}.

For $T\in (0,\infty)$, 
we consider the one-dimensional backward Riccati equation
\begin{equation}
 \dot{R}(t)-a_1(t)R^2(t)+2a_2(t)R(t)+a_3(t)=0, 
 \quad 0\le t\le T,
 \quad R(T)=0,
 \label{eq:B1}
\end{equation}
where
\begin{align}
&\mbox{$a_i(\cdot)\in C([0,\infty)\to \mathbf{R})$ for $i=1,2,3$,}
\label{eq:B2}\\
&\mbox{$a_1(t)\ge 0$ for $t\ge 0$,}
\label{eq:B3}\\
&\mbox{for $i=1,2,3$, $a_i(t)$ converges to $\bar{a}_i$ 
exponentially fast as $t\to\infty$,}
\label{eq:B4}\\
&\mbox{$\bar{a}_1>0$ and $\bar{a_2}^2+\bar{a}_1\bar{a}_3>0$.}
\label{eq:B5}
\end{align}
By (\ref{eq:B5}), we may write $\bar{R}$ for the larger solution to 
the quadratic equation
\[
\bar{a}_1\bar{R}^2-2\bar{a}_2\bar{R}-\bar{a}_3=0.
\]
Recall $\Delta$ from (\ref{eq:2.17}).

\begin{thm}[Nagai and Peng \cite{NP}, Section 5]\label{thm:B1}
We further assume
\begin{equation}
\mbox{$a_3(t)\ge 0$ for $t\ge 0$.}
\label{eq:B6}
\end{equation}
Then, for $T\in (0,\infty)$, {\rm (\ref{eq:B1})} has a unique 
nonnegative solution $R(t)\equiv R(t;T)$, and it satisfies the following:
\begin{enumerate}
\item $R(t;T)$ is bounded in $\Delta$.
\item $\lim_{T-t\to\infty,\ t\to\infty}R(t;T)=\bar{R}$.
\item For $\delta, \epsilon\in (0,\infty)$ such that $\delta+\epsilon<1$, 
\[
\lim_{T\to\infty}\sup_{\delta T\le t\le (1-\epsilon)T}|R(t;T)-\bar{R}|=0.
\]
\end{enumerate}
\end{thm}

When the condition (\ref{eq:B6}) is lacking, we have the following:

\begin{thm}\label{thm:B2}
We assume that, for every $T\in (0,\infty)$, the equation 
{\rm (\ref{eq:B1})} has a solution $R(t)\equiv R(t;T)$ 
that is bounded in $\Delta$. 
Then {\rm (i)--(iii)} in Theorem {\rm \ref{thm:B1}} hold.
\end{thm}

The proof of Theorem \ref{thm:B2} is almost the same as that of Theorem 
\ref{thm:B1} in \cite{NP}, whence we omit it.

We turn to the one-dimensional backward linear differential 
equation
\begin{equation}
\label{eq:B7}
 \dot{v}(t)-b_1(t;T)v(t)+b_2(t;T)=0, 
 \quad 0\le t\le T,
 \quad v(T)=0,
\end{equation}
where
\begin{gather}
\mbox{$b_i(\cdot;T)\in C([0,T]\to \mathbf{R})$ for $T\in (0,\infty)$ and 
$i=1,2$,}
\label{eq:B8}\\
\lim_{T-t\to\infty,\ t\to\infty}b_i(t;T)=\bar{b}_i\ \ \mbox{for}\ \ i=1,2,
\label{eq:B9}\\
\bar{b}_1>0.
\label{eq:B10}
\end{gather}

\begin{thm}[\cite{NP}, Section 5]\label{thm:B3}
For $T\in (0,\infty)$, write $v(t)\equiv v(t;T)$ for the solution to 
{\rm (\ref{eq:B7})}. Let $\bar{v}$ be the solution of the linear equation 
$\bar{b}_1\bar{v}-\bar{b}_2=0$. Then
\begin{enumerate}
\item $v(t;T)$ is bounded in $\Delta$.
\item $\lim_{T-t\to\infty,\ t\to\infty}v(t;T)=\bar{v}$.
\item For $\delta, \epsilon\in (0,\infty)$ such that $\delta+\epsilon<1$,
\[
\lim_{T\to\infty}\sup_{\delta T\le t\le (1-\epsilon)T}|v(t;T)-\bar{v}|=0.
\]
\end{enumerate}
\end{thm}


\section{Parameter estimation}

In this appendix, we use the special case of our model $\mathcal{M}$ 
in which $\sigma_{ij}(t)$'s are constants, i.e.,
\[
\sigma_{ij}(t)=\sigma_{ij},\quad t\ge 0,\quad i,j=1,\dots,n.
\]
We explain how we can statistically estimate 
the parameters $\sigma_{ij}$, $p_i$ and $q_i$ 
from stock price data. 
This problem, for the univariate case $n=1$, is discussed in 
\cite{AIP, INA}. Here we are interested in the multivariate case $n\ge 2$. 
As for the expected rates of return $\mu_i$, 
there is as usual a structural difficulty 
in the statistical estimation of them 
(cf.\ Luenberger \cite[Chapter 8]{L}), 
whence we do not discuss it here.

From (\ref{eq:1.5}), we see that
\begin{equation}
E[Y_j^2(t)]/t=f(t;p_j,q_j),\quad t>0,\quad j=1,\dots,n,
\label{eq:C1}
\end{equation}
where
\[
f(t;p,q):=\frac{q^{2}}{(p+q)^{2}}+\frac{p(2q+p)}{(p+q)^{3}}
\cdot \frac{(1-e^{-(p+q)t})}{t},\quad t>0
\]
(cf.\ \cite{AI}, Examples 4.3 and 4.5). Notice that $f(t;0,q)=1$. 
From (\ref{eq:1.6}) or (\ref{eq:1.7}) 
and the It\^o formula, the solution $S(t)=(S_1(t),\dots,S_n(t))'$ to 
(\ref{eq:1.1}) is given by
\begin{equation}
S_i(t)=s_i\exp\left[\sum\nolimits_{j=1}^n\sigma_{ij} Y_j(t)
+\int_{0}^{t}\left\{\mu_i(s)-\frac{1}{2}\sum\nolimits_{j=1}^n
\sigma_{ij}^2\right\}ds\right],
\quad t\ge 0,
\label{eq:C2}
\end{equation}
for $i=1,\dots,n$. Since $Y(t)$ has stationary increments, we may define
\[
V_{ij}(t-s)
:=\frac{1}{t-s}\mathrm{cov}\left\{\log\frac{S_i(t)}{S_i(s)}, 
\log \frac{S_j(t)}{S_j(s)}\right\},\quad t>s\ge 0, 
\quad i,j=1,\dots,n,
\]
where cov$(\cdot, \cdot)$ denotes the covariance with respect to the 
physical probability measure $P$. 
By (\ref{eq:C1}) and (\ref{eq:C2}), we see that
\begin{gather*}
V_{ij}(t)
=\sum\nolimits_{m=1}^{n}\sigma_{im}\sigma_{jm}f(t;p_m,q_m),
\quad t>0,\quad i,j=1,\dots,n.
\end{gather*}

Suppose that we are given data consisting of closing prices of $n$ assets 
observed at 
a time interval of $N$ consecutive trading days. 
For $m=1,\dots,N$ and $i=1,\dots,n$, 
we denote by $s_{i}(m)$ the price of the $i$th asset 
on the $m$th day. Notice that here the time unit is the day. 
Pick $M<N$, and, for $t\in\{1,\dots,M\}$, 
define $u_i(m)\equiv u_{i}(m,t)$ by 
\[
u_{i}(m):=\log\frac{s_i(m+t)}{s_{i}(m)}, \quad m=1,2,\dots ,N-t.
\]
For $t=1,\dots,M$, we consider the estimator
\begin{equation}
v_{ij}(t):=100^2\frac{252}{t(N-t-1)}
\sum\nolimits_{m=1}^{N-t}
\{u_{i}(m)-\overline{u}_i\}\{u_{j}(m)-\overline{u}_j\}
\label{eq:C3}
\end{equation}
of $V_{ij}(t)$, 
where $\overline{u}_i:=(N-t)^{-1}\sum_{m=1}^{N-t}u_{i}(m)$. 
The number 252, which is the average number of trading days in one 
year, 
converts the return into that per annum, while the number 100 
gives the return in percentage.

We estimate the values of the parameters 
$\sigma_{ij}$, $p_i$ and $q_i$ by nonlinear least squares. 
More precisely, we search for the values of them such that the following 
least squares error is minimized:
\[
\sum_{t=1}^M\sum_{i=1}^n\sum_{j=1}^n\left[V_{ij}(t)-v_{ij}(t)\right]^2.
\]

\begin{figure}[hbtp]
\begin{center}
\includegraphics
[width=350pt,height=240pt]
{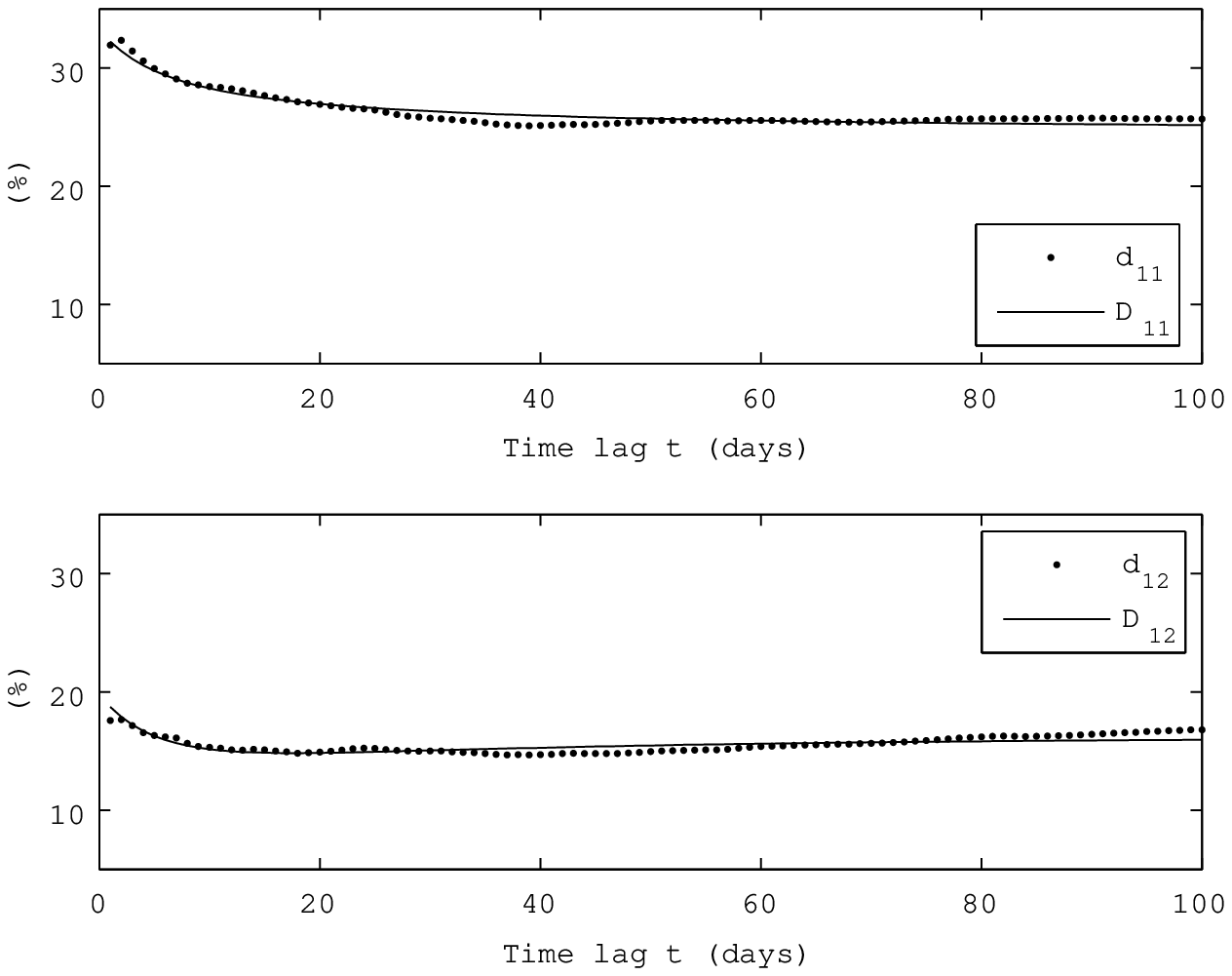}
\end{center}
\caption{$d_{11}(t)$ vs.\ fitted $D_{11}(t)$ and $d_{12}(t)$ vs.\ fitted 
$D_{12}(t)$.}
\label{fig:C1}
\end{figure}

\begin{figure}[hbtp]
\begin{center}
\includegraphics
[width=350pt,height=240pt]
{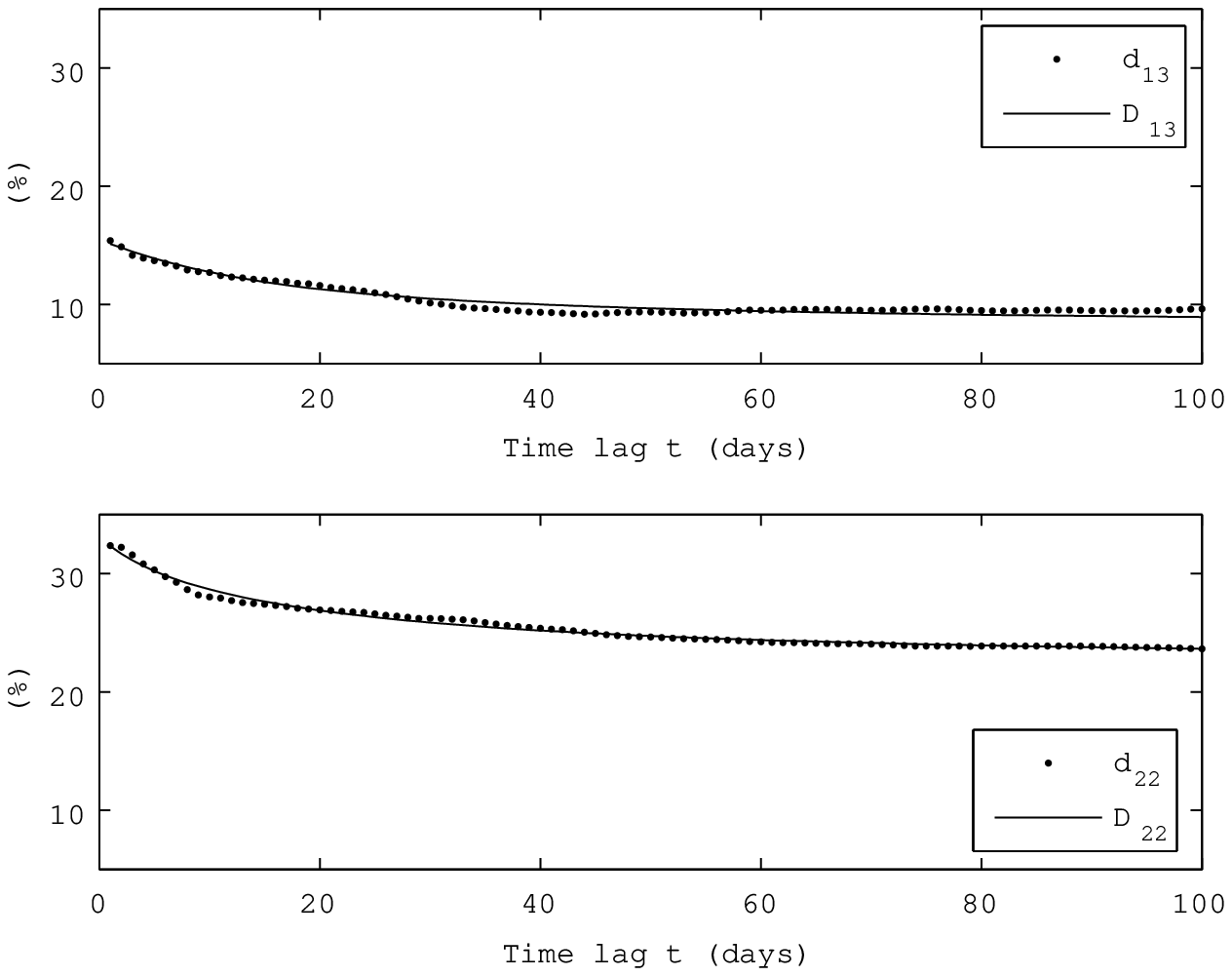}
\end{center}
\caption{$d_{13}(t)$ vs.\ fitted $D_{13}(t)$ and $d_{22}(t)$ vs.\ fitted 
$D_{22}(t)$.}
\label{fig:C2}
\end{figure}

\begin{figure}[hbtp]
\begin{center}
\includegraphics
[width=350pt,height=240pt]
{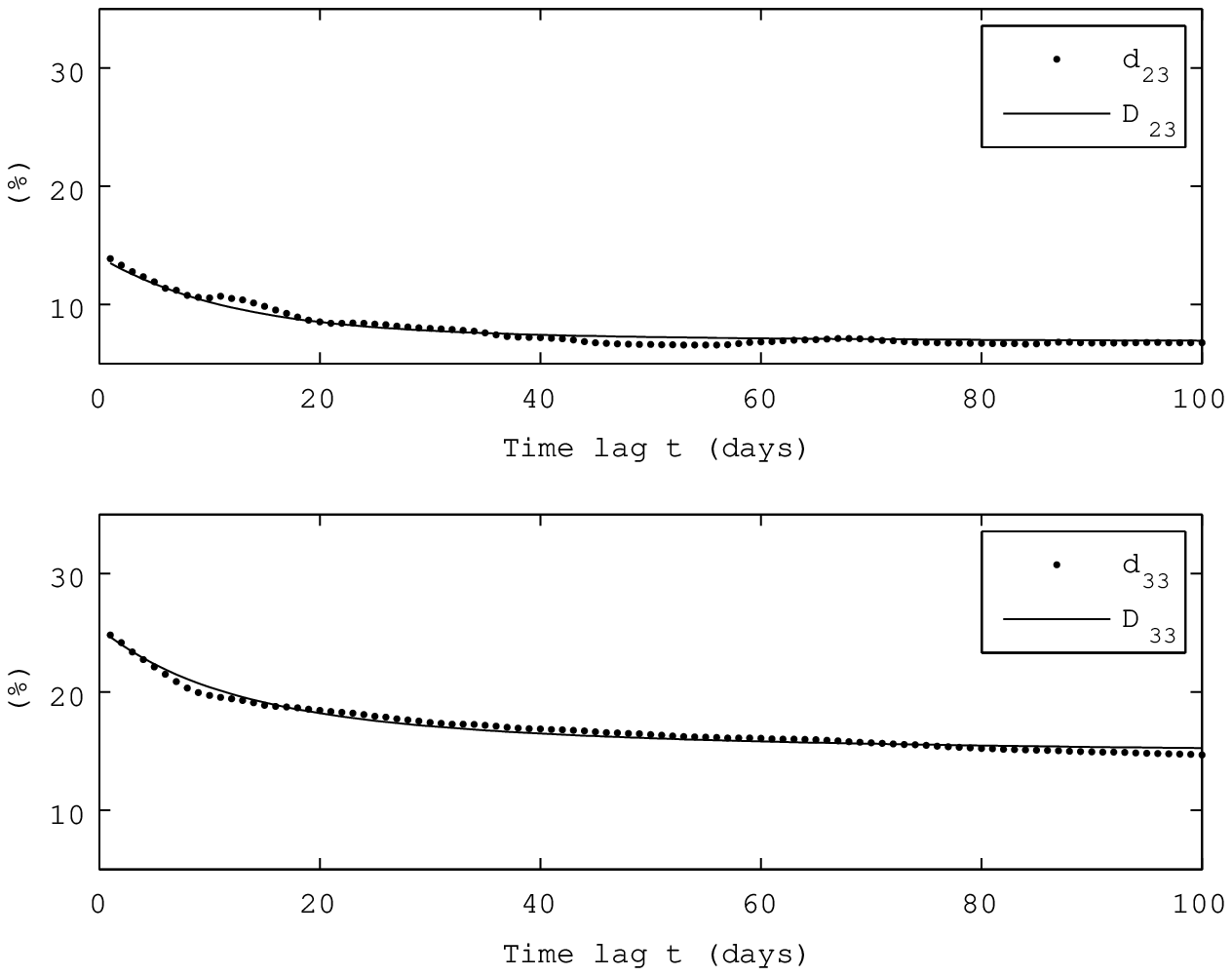}
\end{center}
\caption{$d_{23}(t)$ vs.\ fitted $D_{23}(t)$ and $d_{33}(t)$ vs.\ fitted 
$D_{33}(t)$.}
\label{fig:C3}
\end{figure}

We show numerical results obtained from the following daily stock 
prices from September 18, 1995, through September 16, 2005:
\[
S_1: \mbox{Pfizer Inc.},\quad S_2: \mbox{Wal-Mart Stores Inc.},\quad 
S_3: \mbox{Exxon Mobil Corp.}
\]
Here we use closing prices adjusted for dividends and splits, which are 
available at Yahoo! Finance \cite{Yahoo}, 
rather than actually observed closing prices. 
In this example, we have
\[
n=3,\quad N=2519,\quad M=100.
\]
The estimated values of the parameters are as follows:
\begin{equation*}
\left[
\begin{matrix}
\sigma_{11} & \sigma_{12} & \sigma_{13}\cr
\sigma_{21} & \sigma_{22} & \sigma_{23}\cr
\sigma_{31} & \sigma_{32} & \sigma_{33}
\end{matrix}
\right]
=
\left[
\begin{matrix}
28.7 & -14.1 & 9.1\cr
20.4 &  22.3 & 13.4\cr
-1.8 & -4.6  & 24.9
\end{matrix}
\right],
\quad 
\left[
\begin{matrix}
p_1 & q_1 \cr
p_2 & q_2 \cr
p_3 & q_3 
\end{matrix}
\right]
=
\left[
\begin{matrix}
0.086 & 0.305 \cr
0.261 & 0.044 \cr
0.076 & 0.098
\end{matrix}
\right].
\end{equation*}
Using the signed square root 
$\mathrm{SSR}[x]:=\mathrm{sign}(x)\sqrt{\vert x\vert}$, 
we write
\[
D_{ij}(t):=\mathrm{SSR}[V_{ij}(t)],\quad 
d_{ij}(t):=\mathrm{SSR}[v_{ij}(t)],\quad 
t>0,\quad i,j=1,\dots,n.
\]
In Figures \ref{fig:C1}--\ref{fig:C3}, the dotted lines are 
the graphs of $d_{ij}(t)$'s, while the corresponding solid lines represent 
those of $D_{ij}(t)$'s that are obtained by using the nonlinear least 
squares above. We see that the fitted functions 
$D_{ij}(t)$ simultaneously approximate the corresponding sample values 
$d_{ij}(t)$ well for this data set. 
We have repeated this procedure for various data sets and 
obtained reasonably good fits in most cases.

\ 

{\bf Acknowledgements.}\quad 
We express our gratitude to Jun Sekine who suggested the use of a 
Cameron--Martin type formula in the proof of Theorem \ref{thm:2.3}. 
The work of the second author was partially supported by a Research 
Fellowship of the Japan Society for the Promotion of Science for 
Young Scientists.

\

\end{document}